\documentclass[a4paper,11pt]{article}
\usepackage{latexsym,amsfonts,amsmath}
\usepackage{amssymb}
\usepackage{color,soul} 
\usepackage{xcolor}
\usepackage[mathscr]{eucal}
\usepackage{enumerate}
\usepackage{enumitem}

\usepackage{pgf}
\usepackage{tikz}
\usetikzlibrary{arrows,automata}

\usepackage{fullpage}

\definecolor{Dgreen}{rgb}{102,255,102}
\definecolor{Purp}{rgb}{102,0,204}

\def\tto{\twoheadrightarrow}

\def\1{\mathbf{1}}
\def\0{\mathbf{0}}

\def\NN{\mathbb{N}}

\def\EE{\mathbb{E}}
\def\RR{\mathbb{R}}

\def\car{{\mathcal{C}ar}}

\def\XX{\mathbf{X}}
\def\AA{\mathbf{A}}
\def\KK{\mathbf{K}}
\def\SS{\mathbf{S}}

\newcommand{\widebar}[1]{\overline{#1}}

\newcommand{\bcal}[1]{\boldsymbol{\mathcal{#1}}}

\newcommand{\bfrak}[1]{\boldsymbol{\mathfrak{#1}}}

\def\argmax{\mathop{\hbox{\rm arg max}}\limits}

\def\liminf{\mathop{\underline{\lim}}}

\def\ds{\displaystyle}

\def\wstar{\ {\stackrel{*}{\rightharpoonup}\ }}

\definecolor{Dgreen}{rgb}{0,0.5,0}

\newcounter{hypot}

    \newenvironment{hypot}{\begin{list}
      {\hspace{\labelsep} \bfseries Assumption $\boldsymbol{\mathrm{\Alph{hypot}}}$ }
      {\leftmargin=0pt
       \labelwidth=0cm
       \refstepcounter{hypot}
       \def\makelabel##1{##1}}}{\end{list}}    

\newcounter{assump}
\renewcommand{\theassump}{$(\mathrm{\Alph{hypot}}_\arabic{assump})$}

\newenvironment{assump}{\begin{list}
      {\hspace{\labelsep} \theassump}
      {\leftmargin=0pt
       \labelwidth=0cm
       \usecounter{assump}
       }}{\end{list}}
       
\newcounter{hypotcal}

    \newenvironment{hypotcal}{\begin{list}
      {\hspace{\labelsep}\bfseries Assumption $\bcal{\Alph{hypotcal}}$}
      {\leftmargin=0pt
       \labelwidth=0cm
       \refstepcounter{hypotcal}
       \def\makelabel##1{##1}}}{\end{list}}    

\newcounter{assumpcal}
\renewcommand{\theassumpcal}{$(\mathcal{\Alph{hypotcal}}_\arabic{assumpcal})$}

\newenvironment{assumpcal}{\begin{list}
      {\hspace{\labelsep} \theassumpcal}
      {\leftmargin=0pt
       \labelwidth=0cm
       \usecounter{assumpcal}
       }}{\end{list}}
       
\newcounter{hypotprim}

    \newenvironment{hypotprim}{\begin{list}
      {\hspace{\labelsep}\bfseries Assumption $\boldsymbol{\mathrm{\Alph{hypotprim}}'}$}
      {\leftmargin=0pt
       \labelwidth=0cm
       \refstepcounter{hypotprim}
       \def\makelabel##1{##1}}}{\end{list}}    

\newcounter{assumpprim}
\renewcommand{\theassumpprim}{$(\mathrm{\Alph{hypotprim}}'_\arabic{assumpprim})$}

\newenvironment{assumpprim}{\begin{list}
      {\hspace{\labelsep} \theassumpprim}
      {\leftmargin=0pt
       \labelwidth=0cm
       \usecounter{assumpprim}
       }}{\end{list}}

\begin{document}
\newtheorem{theorem}{Theorem}[section]
\newtheorem{proposition}[theorem]{Proposition}
\newtheorem{lemma}[theorem]{Lemma}
\newtheorem{corollary}[theorem]{Corollary}
\newtheorem{definition}[theorem]{Definition}
\newtheorem{remark}[theorem]{Remark}
\newtheorem{conjecture}[theorem]{Conjecture}
\newtheorem{assumption}[theorem]{Assumption}

\title{Nash equilibria for total expected reward  absorbing Markov games: the constrained and unconstrained cases\footnote{Supported by grant PID2021-122442NB-I00 from the Spanish \textit{Ministerio de Ciencia e Innovaci\'on.}}}

\author{Fran\c{c}ois Dufour\footnote{
Institut Polytechnique de Bordeaux; INRIA Bordeaux Sud Ouest, Team: ASTRAL; IMB, Institut de Math\'ematiques de Bordeaux, Universit\'e de Bordeaux, France \tt{francois.dufour@math.u-bordeaux.fr}}
\and
Tom\'as Prieto-Rumeau\footnote{Statistics Department, UNED, Madrid, Spain. e-mail: {\tt{tprieto@ccia.uned.es}}\quad {(Author for correspondence)}}}
\maketitle
\begin{abstract}
We consider a nonzero-sum $N$-player  Markov  game on an abstract measurable state space with compact metric action spaces. The payoff functions are bounded Carath\'eodory functions and the transitions of the system are assumed to have a density function satisfying some continuity conditions. The optimality criterion of the players is given by a total expected payoff on an infinite discrete-time horizon. Under the condition that the game model is absorbing, we establish the existence of Markov strategies that are a noncooperative equilibrium in the family of all history-dependent strategies of the players for both the constrained and the unconstrained problems, 
We obtain, as a particular case of results, the existence of Nash equilibria for discounted constrained and unconstrained game models.
\end{abstract}
{\small 
}

\section{Introduction}\label{sec-1}
The topic of  noncooperative  games has been extensively studied in the last decades and research on this subject has spread in many directions. Here, we are interested in nonzero-sum Markov games; namely, we deal with a stochastic dynamic game on an infinite time horizon: the state of the system evolves according to a stochastic Markov-like kernel, players take their actions  after each transition of the system, and their goal is to maximize a given optimality criterion in a noncooperative way. The primary goal is to establish  the existence of a Nash equilibrium: that is, a strategy for each of the players in such a way that none of the players can improve his payoff by using another strategy. The interested reader can consult, for instance, the survey \cite{jaskiewicz18} to have an overview of this topic. 

When dealing with the discounted payoff optimality criterion, the usual technique is  the dynamic programming approach. This consists in considering a one-shot game (with a single decision epoch) and then establish the existence of a fixed point for some suitably defined selectors, which will yield the Nash equilibrium. Such results have been obtained under various hypotheses in, e.g., \cite{nowak03,nowak92}, and the most refined and elegant conditions  have been proposed in \cite{he-sun17} by introducing  so-called ``decomposable coarser transition kernels''.

A natural generalization of the game models described so far is to consider games with \textit{constraints}. In this case, the players try to maximize their payoff function subject to the condition that some constraints (related to some other payoff functions) must be satisfied. A noncooperative equilibrium is then defined as a set of strategies of the players which satisfy simultaneously all the constraints and for which, in addition, no player can improve his payoff when unilaterally modifying his strategy while still satisfying his own constraints. Such constrained game models were studied first in 
 \cite{Altman-Shwartz-96} for a model with  finite state and actions spaces, and then generalized in \cite{Mena-OHL-constrained-games} for game models with countable state space and compact action spaces, under some conditions which make the countable state space model nearly finite  in each transition of the system.
 It is not possible, in general, to solve these constrained game models by using the dynamic programming approach. The linear programming approach (which consists in considering the spaces of occupation measures associated to the strategies of the players) appears to be well suited to study constrained problems.

Extending this linear programming approach, in the context of obtaining Nash equilibria for game models, from the finite or countable state space cases to a more general state space entails, however, serious technical difficulties. 
This was achieved in \cite{dufour22} and, based on this reference, 
there has been a recent growing interest in the study of constrained games. In \cite{dufour22}, a game model on a measurable state space with ARAT (additive reward, additive transition) structure under the discounted optimality criterion was studied. The approach consists in  defining a correspondence on a space of measures, endowed with the weak-strong topology, which is shown to have a fixed point. From this fixed point, optimal constrained strategies of the players are obtained.  The approach in~\cite{dufour22} combines also the use of Young measures to identify Markov strategies of the players. With this approach, in the reference \cite{Jaskiewicz-Nowak-AMO-22}, a countable state space discounted game model is studied, dropping the ARAT condition. Also, in \cite{Jaskiewicz-Nowak-AMO-23}, the authors study a discounted constrained game on a general state space and consider the weaker notion of an approximate equilibria. 

The above cited references are, generally, concerned with  games under the \textit{discounted payoff} optimality criterion.
In this paper, we shall consider games under the total expected payoff optimality criterion,  for constrained and unconstrained games, both of which have
received much smaller attention than the discounted payoff counterpart.
Indeed,
as far as we know, only a few references deal with that topic: in \cite{Padilla-total-game}, the existence of an $\epsilon$-equilibrium for a countable state space game is established, while in the references \cite{cavazos-total-game-MMOR,cavazos-total-game-Kyb}, a stopping zero-sum game (in which one player is allowed to stop the evolution of the system) with countable state space under the total expected payoff criterion  is studied, and in \cite{Cortes-Kyb} for a finite state game. Summarizing, there do not exist general existence results for Nash equilibria for the total expected payoff criterion, even in the case of a countable state space.

The extension from the  discounted optimality criterion to the total expected payoff optimality criterion is far from being straightforward. First of all, the dynamic programming approach using the one-shot game does not work in our context: indeed, the usual technique ---see, e.g., \cite{nowak92}--- consists in showing that the Bellman operator maps the final payoff function, which typically belongs to a compact subset of the space of functions $L^\infty$, into itself when computing the one-shot optimal payoff. To use this technique, however, the discount factor plays a crucial role and, for the total expected payoff criterion, the dynamic programming approach does not yield an operator mapping a compact set of functions into itself. 
Secondly, while, under the discounted optimality criterion, the occupation measures are probability or finite measures, in the context of the total expected payoff criterion, the occupation measures may be infinite. Ensuring finiteness of these occupation measures and establishing compactness properties of these spaces of measures becomes more technically demanding  and some  additional conditions must be imposed. More precisely, we will need to assume that the game model is uniformly absorbing, which means that the dynamic system enters into some subset $\Delta$ of the state space (in which no further reward is earned) in a uniformly bounded expected time, with the queues of the hitting time of $\Delta$ converging uniformly to zero as well (a precise definition will be given in the text).  As already mentioned, our techniques herein will allow to deal, at the same time, with constrained and unconstrained game models. 

In this paper we will consider  an $N$-player Markov  Markov game under the total expected reward optimality criterion. The state space is an abstract measurable space and the action spaces are compact metric spaces. The payoff functions are assumed to be bounded Carath\'eodory functions and we impose that the transition probabilities have a density function with respect to some reference probability measure, and that they satisfy some suitable $L^1$-continuity properties. Some uniformly absorbing requirement must be imposed to ensure finiteness of the payoff criteria and compactness of the occupation measures, plus the usual Slater condition. It is worth stressing that we do not need the ARAT separation property and, instead, we impose some minimal sufficient conditions ensuring the continuity of the transition and cost functionals. 

The rest of the paper is organized as follows. In the remaining of this section we introduce some notation and recall some standard results that will be useful in the sequel. Section \ref{sec-2} is devoted to define the constrained and unconstrained game models, and  to propose some basic assumptions. Our main results in the paper are stated in Section \ref{sec-3}. 
In Section \ref{section-Assumptions-YoungMeasures} we study the occupation measures and introduce the spaces of Young measures, which shall be identified with Markov stationary strategies of the players. Some useful continuity results relating narrow convergence of Young measures and weak-strong convergence of measures are established. Our main results on the existence of constrained and unconstrained equilibria are proved in Section \ref{section4}.

\paragraph{Notation and terminology.}
A metric space $\SS$ will be always endowed with its Borel $\sigma$-algebra  $\bfrak{B}(\SS)$. On the product of a finite number of metric spaces $\SS=\SS^1\times\ldots\times\SS^N$, we will consider the product topology which makes the product again a metric space. If the metric spaces $\SS^i$ are separable, then we have $\bfrak{B}(\SS)=\bfrak{B}(\SS^1)\otimes\ldots\otimes\bfrak{B}(\SS^N)$.

On a measurable space $(\mathbf{\Omega},\mathcal{F})$ we will consider the set of finite signed measures $\bcal{M}(\mathbf{\Omega})$, the set of finite nonnegative measures
$\bcal{M}^+(\mathbf{\Omega})$, and the set of probability measures~$\bcal{P}(\mathbf{\Omega})$.
For a set $\Gamma\in\mathcal{F}$, we denote by $\mathbf{I}_{\Gamma}:\Omega\rightarrow\{0,1\}$ the indicator function of the set $\Gamma$, that is,
$\mathbf{I}_{\Gamma}(\omega)=1$ if and only if $\omega\in\Gamma$.
For $\omega\in\mathbf{\Omega}$, we write $\delta_{\{\omega\}}$ for the Dirac probability measure at $\omega$ defined on $(\mathbf{\Omega},\mathcal{F})$ by
$\delta_{\{\omega\}}(B)=\mathbf{I}_{B}(\omega)$ for any $B\in\mathcal{F}$. If $\mu\in\bcal{M}(\mathbf{\Omega})$ and $\Gamma\in\mathcal{F}$, we denote by
$\mu_{\Gamma}$ the measure on $(\mathbf{\Omega},\mathcal{F})$ defined by $\mu_{\Gamma}(B)=\mu(\Gamma\cap B)$ for $B\in\mathcal{F}$.
The trace $\sigma$-algebra of a set $\Gamma\subseteq\mathbf{\Omega}$ is denoted by~$\mathcal{F}_{\Gamma}$.
On $\bcal{P}(\mathbf{\Omega})$, the $s$-topology is the coarsest topology that makes  $\mu\mapsto \mu(D)$ continuous for every $D\in\mathcal{F}$.

Given a measurable space  $(\mathbf{\Omega},\mathcal{F})$ and $\lambda\in\bcal{P}(\mathbf{\Omega})$, we will 
denote by $L^{1}(\mathbf{\Omega},\mathcal{F},\lambda)$  the family of measurable functions  (identifying those which are $\lambda$-a.s. equal)  $f:\mathbf{\Omega}\rightarrow\RR$ which are $\lambda$-integrable, i.e., $\|f\|_1=\int_\XX |f(x)|\lambda(dx)<\infty$. 
Also, let  $L^{\infty}(\mathbf{\Omega},\mathcal{F},\lambda)$ be the set of $\lambda$-essentially bounded measurable functions $f:\mathbf{\Omega}\rightarrow\RR$  (again, we identify functions that  coincide $\lambda$-a.s.). We will denote by  $\|f\|_\infty$ the corresponding  essential supremum. On $L^{\infty}(\mathbf{\Omega},\mathcal{F},\lambda)$ we will consider the  weak$^{*}$ topology, that is, we have $f_n\wstar f$ whenever
$$\int_\mathbf{\Omega} f_nh \,d\lambda\rightarrow \int_\mathbf{\Omega} fh\,d\lambda\quad\hbox{for every $h\in L^{1}(\mathbf{\Omega},\mathcal{F},\lambda)$}.$$

Let $(\mathbf{\Omega},\mathcal{F})$ and $(\widetilde{\mathbf{\Omega}},\widetilde{\mathcal{F}})$ be two measurable spaces.
A kernel on $\widetilde{\mathbf{\Omega}}$ given $\mathbf{\Omega}$ is a mapping
$Q:\mathbf{\Omega}\times\widetilde{\mathcal{F}}\rightarrow\RR^+$ such that $\omega\mapsto Q(B|\omega)$ is 
measurable on $(\mathbf{\Omega},\mathcal{F})$ for every $B\in\widetilde{\mathcal{F}}$,   and  $B\mapsto Q(B|\omega)$ is in $\bcal{M}^+(\widetilde{\mathbf{\Omega}})$
for every $\omega\in\mathbf{\Omega}$. If $Q(\widetilde{\mathbf{\Omega}}|\omega)=1$ for all $\omega\in\mathbf{\Omega}$ then we say that $Q$ is a \textit{stochastic kernel}.
We write $\mathbb{I}_{\Gamma}$ for the kernel on $\mathbf{\Omega}$ given $\mathbf{\Omega}$ defined by
$\mathbb{I}_{\Gamma}(B|\omega)=\mathbf{I}_{\Gamma}(\omega) \delta_{\{\omega\}}(B)$ for $\omega\in\mathbf{\Omega}$ and $B\in\mathcal{F}$.
Let $Q$ be a stochastic kernel on $\widetilde{\mathbf{\Omega}}$ given $\mathbf{\Omega}$.
For  a bounded measurable function $f:\widetilde{\mathbf{\Omega}}\rightarrow\RR$, we will denote by $Qf:\mathbf{\Omega}\rightarrow\RR$ the measurable function
$$Qf(\omega)=\int_\mathbf{\Omega'} f(\widetilde{\omega})Q(d\widetilde{\omega}|\omega)\quad\hbox{for $\omega\in\mathbf{\Omega}$}.$$
For a measure $\mu\in\bcal{M}^{+}(\mathbf{\Omega})$, we denote by $\mu Q$ the finite measure  on $(\widetilde{\mathbf{\Omega}},\widetilde{\mathcal{F}})$ given by
$$B\mapsto \mu Q\,(B)= \int_{\mathbf{\Omega}} Q(B|\omega) \mu(d\omega)\quad\hbox{for $B\in\widetilde{\mathcal{F}}$}.$$
The product of the $\sigma$-algebras $\mathcal{F}$ and $\widetilde{\mathcal{F}}$ is denoted by $\mathcal{F}\otimes\widetilde{\mathcal{F}}$ and consists of the $\sigma$-algebra
generated by the measurable rectangles, that is, the sets of the form $\Gamma\times\widetilde{\Gamma}$ for $\Gamma\in\mathcal{F}$ and
$\widetilde{\Gamma}\in\widetilde{\mathcal{F}}$.
We denote by $\mu\otimes Q$ the unique probability measure (or finite measure) on the product space $(\mathbf{\Omega}\times\widetilde{\mathbf{\Omega}},\mathcal{F}\otimes\widetilde{\mathcal{F}})$ satisfying 
$$(\mu\otimes Q)( \Gamma\times\widetilde{\Gamma})= \int_{\Gamma} Q(\widetilde{\Gamma}|\omega)\mu(d\omega)\quad\hbox{for $\Gamma\in\mathcal{F}$ and
$\widetilde{\Gamma}\in\widetilde{\mathcal{F}}$,}$$
see Proposition III-2-1 in \cite{neveu70} for a proof of existence and uniqueness of such measure.
Let $(\widebar{\mathbf{\Omega}},\widebar{\mathcal{F}})$ be a third measurable space and $R$ a stochastic kernel on $\widebar{\mathbf{\Omega}}$
given $\widetilde{\mathbf{\Omega}}$. Then we will denote by $QR$ the stochastic kernel on $\widebar{\mathbf{\Omega}}$ given $\mathbf{\Omega}$ given by
$$QR(\Gamma |\omega)= \int_{\widetilde{\mathbf{\Omega}}} R(\Gamma | \tilde{\omega}) Q(d\tilde{\omega} | \omega)
\quad\hbox{for $\Gamma\in\widebar{\mathcal{F}}$ and $\omega\in\mathcal{F}$}.$$
Given $\mu\in\bcal{M}(\mathbf{\Omega}\times\widetilde{\mathbf{\Omega}})$,  the  marginal measures are $\mu^{\mathbf{\Omega}}\in\bcal{M}(\mathbf{\Omega})$ and
$\mu^{\widetilde{\mathbf{\Omega}}}\in\bcal{M}(\widetilde{\mathbf{\Omega}})$ defined by 
$\mu^{\mathbf{\Omega}}(\cdot)=\mu(\cdot\times \widetilde{\mathbf{\Omega}})$ and
$\mu^{\widetilde{\mathbf{\Omega}}}(\cdot)=\mu(\mathbf{\Omega}\times\cdot)$. If $\pi$ is a kernel on $\widetilde{\mathbf{\Omega}},\times\widebar{\mathbf{\Omega}}$ given $\mathbf{\Omega}$ the  marginal kernels are $\pi^{\widetilde{\mathbf{\Omega}}}$
and $\pi^{\widebar{\mathbf{\Omega}}}$, respectively defined by 
$\pi^{\widetilde{\mathbf{\Omega}}}=\pi(\cdot\times \widebar{\mathbf{\Omega}}|\omega)$ and
$\pi^{\widebar{\mathbf{\Omega}}}=\pi(\widetilde{\mathbf{\Omega}}\times\cdot |\omega)$ for $\omega\in\mathbf{\Omega}$.

We say that $f:\mathbf{\Omega}\times\SS\rightarrow\SS'$, where~$\SS'$ is a metric space, is a \textit{Carath\'eodory function} if $f(\cdot,s)$ is measurable on~$\mathbf{\Omega}$ for every $s\in\SS$ and $f(\omega,\cdot)$ is continuous on $\SS$ for every $\omega\in \mathbf{\Omega}$. The family of the so-defined Carath\'eodory functions is denoted by $\car(\mathbf{\Omega}\times\SS,\SS')$. The family of Carath\'eodory functions which, in addition, are bounded is denoted by $\car_b(\mathbf{\Omega}\times\SS,\SS')$. 
When the metric space~$\SS$ is separable then any $f\in\car(\mathbf{\Omega}\times\SS,\SS')$ is a jointly measurable function on $(\mathbf{\Omega}\times\mathbf{S},\mathcal{F}\otimes\bfrak{B}(\mathbf{S}))$; see \cite[Lemma 4.51]{aliprantis06}. 

Given $\lambda\in\bcal{P}(\mathbf{\Omega})$, let $\bcal{P}_{\lambda}(\mathbf{\Omega})=\{\eta\in  \bcal{P}(\mathbf{\Omega}): \eta\ll\lambda\}$ be the family of probability  of probability measures which are absolutely continuous with respect to $\lambda$.

If $\SS$ is a Polish space (a complete and separable metric space), on $\bcal{M}(\mathbf{\Omega}\times\SS)$ we will consider the  $ws$-topology (weak-strong topology) which is the coarsest topology for which the mappings
$$
\mu\mapsto \int_{\mathbf{\Omega}\times\SS} f(\omega,s)\mu(d\omega,ds)
$$
for $f\in\car_b(\mathbf{\Omega}\times\SS,\RR)$ are continuous.
There are other equivalent definitions of this topology as discussed, for instance, in \cite[Section 3.3]{florescu12}.

Inequality $\ge$ in $\RR^p$ means a componentwise inequality $\ge$, while the inequality $>$ in $\RR^p$ is a componentwise strict inequality $>$. Let $\1\in\RR^p$ be the vector with all components equal to one.

The next disintegration lemma will be useful in the forthcoming (see Theorem 1 in \cite{valadier73}).
\begin{lemma}[Disintegration lemma]
\label{lemma-disintegration} 
Let $(\mathbf{\Omega},\mathcal{F})$ be a measurable space and let $\SS$ be a Polish space. Let $\varphi:\mathbf{\Omega}\tto\SS$ be a weakly measurable correspondence with nonempty closed values, and let~$\mathbf{K}$ be the graph of the correspondence. For every  $\mu\in\bcal{M}^{+}(\mathbf{\Omega}\times\SS)$ 
such that $\mu(\mathbf{K}^{c})=0$
there exists a stochastic kernel $Q$ on $\SS$ given $\mathbf{\Omega}$ such that 
\begin{equation}\label{eq-dudley-product}
\mu= \mu^{\mathbf{\Omega}}\otimes Q 
\end{equation}
and such that
$Q(\varphi(\omega)|\omega)=1$ for each $\omega\in\mathbf{\Omega}$. Moreover, $Q$ is unique $\mu^{\mathbf{\Omega}}$-almost surely, meaning that if $Q$ and $Q'$ are two stochastic kernels that satisfy \eqref{eq-dudley-product} then for all $\omega$ in a set of $\mu^{\mathbf{\Omega}}$-probability one, the probability measures $Q(\cdot|\omega)$ and $Q'(\cdot|\omega)$ coincide.
\end{lemma}

\section{Definition of the game model}
\label{sec-2}

\subsection{Elements of the noncooperative game model}

Next we give the primitive data of our $N$-person game model.

\begin{enumerate}
\item [(a).] The state space is a measurable space $\XX$ endowed with a $\sigma$-algebra $\bfrak{X}$. 
\item [(b).]  The separable metric space $\AA^i$, with $i\in\{1,\ldots,N\}$, stands for the action space of player $i$. Given any $x\in\XX$, the nonempty measurable set $\AA^i(x)\subseteq \AA^i$ is the set of actions  available to player $i$ at state $x$. We will use the notations
$$\AA=\AA^1\times\ldots\times\AA^N\quad\hbox{and}\quad \AA(x)=\AA^1(x)\times\ldots\times\AA^N(x)\subseteq\AA.$$
A typical element of $\AA$ will be written $a=(a^{1},\ldots,a^{N})$.
\item [(c).] Given  $i\in\{1,\ldots,N\}$, the bounded measurable  functions $r^{i}:\XX\times\AA\rightarrow\RR$ and $c^i:\XX\times\AA\rightarrow\RR^p$ stand for the reward  and constraint functions for  player $i$. The components of $c^i$ will be denoted by $c^{i,j}$ for $1\le j\le p$. The corresponding constraint constant is $\rho^i\in\RR^p$. Here, $p\ge1$ is a fixed integer assumed to be the same for all the players.  We write $\rho=(\rho^1,\ldots,\rho^N)\in\RR^{pN}$.
\item [(d).] The transitions of the system are given by a stochastic kernel $Q$ on $\mathbf{X}$ given $\mathbf{X}\times\mathbf{A}$. 
\item [(e).] The initial distribution of the system is the probability measure $\eta\in\bcal{P}(\XX)$.
\end{enumerate}
We will consider the sets  $$\mathbf{H}_0=\XX\quad\hbox{and}\quad\mathbf{H}_t=(\XX\times\AA)^t\times\XX\ \ \hbox{for $t\ge1$}, $$which are the sets of histories of the state-action process up to time $t\ge0$. 
An element of~$\mathbf{H}_t$ is denoted by $h_t=(x_0,a_0,\ldots,x_{t-1},a_{t-1},x_t)$.
\textit{We note that, throughout this paper, sub-indices will usually refer to the time component $t\ge0$, while super-indices will typically denote the players $i\in\{1,\ldots,N\}$.}

For this model, it is assumed that, at time $t\ge0$, the players choose their actions independently of each other conditionally on the history $h_{t}$ of the system; hence, its noncooperative nature.
This game will be denoted by $\mathcal{G}(\eta,\rho)$. We use this notation because, in the sequel, we will need to vary both the initial distribution and the constraint constants.
\begin{definition}\label{def-policies} Fix a player $i\in\{1,\ldots,N\}$.
\begin{itemize}
\item[(i).] A policy for player $i$  is a sequence $\{\pi_t^i\}_{t\ge0}$ of stochastic kernels on $\AA^i$ given $\mathbf{H}_t$ that verify 
$\pi^i_t(\AA^i(x_t)|x_0,a_0,\ldots,x_t)=1$ for every $t\ge0$ and $h_t=(x_0,a_0,\ldots,x_t)\in\mathbf{H}_t$.
Let $\mathbf{\Pi}^i$ be the family of all policies of player $i$. 
\item[(ii).] Let $\mathbf{M}^i$ be the family of stochastic kernels $\pi^i$ on $\AA^i$ given $\XX$ for which $\pi^i(\AA^i(x)|x)=1$ for every $x\in\XX$. The policy 
$\{\pi^i_t\}_{t\ge0}\in\mathbf{\Pi}_i$ is said to be a stationary Markov policy for player $i$ if, for some $\pi^i\in\mathbf{M}^i$, we have
$$\pi^i_t(\cdot|x_0,a_0,\ldots,x_t)=\pi^i(\cdot|x_t)\quad\hbox{for every $t\ge0$ and $h_t=(x_0,a_0,\ldots,x_t)\in\mathbf{H}_t$}.$$
\end{itemize}
\end{definition}
Our Assumption \ref{Assumption-control} below will ensure that these sets of policies  are nonempty.
We will usually refer to $\mathbf{\Pi}^i$ as to  the class of history-dependent policies of player $i$. The family of history-dependent policies for the players is $\mathbf{\Pi}=\mathbf{\Pi}^1\times\ldots\times \mathbf{\Pi}^N$. We will say that $\pi\in\mathbf{\Pi}$ is a strategy profile. 
We can identify the class of stationary Markov policies for player $i$ with $\mathbf{M}^i$ itself, and so we will write $\mathbf{M}^i\subseteq\mathbf{\Pi}^i$. 
Similarly, we introduce the notation $\mathbf{M}=\mathbf{M}^1\times\ldots \times\mathbf{M}^N$ for the class of stationary Markov profiles of the players.  

\paragraph{The $-i$ notation.}

Given $\pi=(\pi^1,\ldots,\pi^N)$ and some $i\in\{1,\ldots,N\}$, let
$$\pi^{-i}=(\pi^1,\ldots,\pi^{i-1},\pi^{i+1},\ldots,\pi^N)\in\mathbf{\Pi}^1\times\ldots\times \mathbf{\Pi}^{i-1}\times \mathbf{\Pi}^{i+1}\times\ldots\times \mathbf{\Pi}^{N}.$$ 
In addition, given $\sigma\in\mathbf{\Pi}^i$, we will use the notation $(\pi^{-i},\sigma)$ to denote the strategy profile in $\mathbf{\Pi}$ for which player $i$ uses the policy $\sigma$ and the remaining players use the policies $\pi^j$ for $j\neq i$.
Similarly, we will use notations such as $\mathbf{\Pi}^{-i}$ and $\mathbf{M}^{-i}$ to consider the product spaces of all the $\mathbf{\Pi}^j$ and $\mathbf{M}^j$ except $\mathbf{\Pi}^{i}$ and $\mathbf{M}^{i}$,  respectively.

\paragraph{Construction of the state-action process.}
The canonical space 
 $\mathbf{H}_\infty=(\XX\times\AA)^\NN$ is endowed with the product $\sigma$-algebra $(\bfrak{X}\otimes\bfrak{B}(\AA))^\NN$. Let  $(X_t,A_t)_{t\ge0}$ be the corresponding coordinates mappings with $A_t=(A_t^1,\ldots,A_t^N)$. We shall use the notation  $H_t=(X_0,A_0,\ldots,X_t)$ for $t\geq1$ and $H_{0}=X_{0}$.
 Let us consider an initial probability measure $\eta\in\bcal{P}(\XX)$ and a strategy profile $\pi\in\mathbf{\Pi}$. There exists a unique probability measure 
$\mathbb{P}_{\eta,\pi}$ on $\mathbf{H}_\infty$ such that for every $B\in\bfrak{X}$, $C^i\in\bfrak{B}(\AA^i)$ for $i=1,\ldots,N$, and $t\ge0$ we have: (i):
$\mathbb{P}_{\eta,\pi}\{X_0\in B\}=\eta(B)$;
\begin{equation}\label{eq-dynamic-game}
\text{(ii): }\mathbb{P}_{\eta,\pi}(A^1_t\in C^1,\ldots,A^N_t\in C^N \mid H_t)=\pi_t^1(C^1|H_t)\cdots\pi_t^N(C^N|H_t),
\end{equation}
and (iii):
$\mathbb{P}_{\eta,\pi}(X_{t+1}\in B| H_t,A_t)=Q(B|X_t,A_t)$.
We denote by  $\mathbb{E}_{\eta,\pi}$
the expectation operator associated to $\mathbb{P}_{\eta,\pi}$. If the initial distribution is the Dirac measure $\delta_{x}$ concentrated at a given state $x\in\XX$ we will simply write $\mathbb{P}_{x,\pi}$ and $\mathbb{E}_{x,\pi}$.

\subsection{Correlated strategies, absorbing models, and Nash equilibria} 
In Definition \ref{def-policies} it is assumed that the players choose their actions independently of each other.
It will be technically useful, however, to introduce \textit{correlated} strategies for which the players can take \textit{dependent actions}.
The notion of correlated strategies plays a very important role in the analysis of this type of game. 
In particular, it will allow to introduce the set of possible answers for each player (see Definition \ref{Definition-Otilde-i}) which are defined from the occupation measures of the process generated precisely from these correlated strategies
(see Definition \ref{def-occupation-measure}).
To show our main results, we will have to show that this set of occupation measures is bounded and compact, which leads to a notion of absorbing model also defined on the basis of correlated strategies.

\begin{definition}\label{def-correlated}
\begin{enumerate}
\item[(i).]
Let  $\tilde{\mathbf{\Pi}}$ be the set of correlated strategies defined as follows: we say that $\pi=\{\pi_t\}_{t\ge0}$ is in $ \tilde{\mathbf{\Pi}}$ if, for every $t\ge0$ and $h_t=(x_0,a_0,\ldots,x_t)\in\mathbf{H}_t$, we have that $\pi_t$ is a stochastic kernel on $\AA$ given $\mathbf{H}_t$ that verifies
$\pi_t(\AA(x_t)|x_0,a_0,\ldots,x_t)=1$.
\item[(ii).]
The class of correlated Markov strategies of the players is $\tilde{\mathbf{M}}$ and it is defined as the set of stochastic kernels $\pi$ on $\AA$ given $\XX$ such that $\pi(\AA(x)|x)=1$ for every $x\in\XX$.  
As usual, we will assume that  $\tilde{\mathbf{M}}\subseteq\tilde{\mathbf{\Pi}}$.
\item[(iii).] Given an initial distribution $\eta\in\bcal{P}(\XX)$ and a correlated strategy $\pi\in\tilde{\mathbf{\Pi}}$, we can construct the state-action process as we did before, except that \eqref{eq-dynamic-game} is replaced with
$\mathbb{P}_{\eta,\pi}(A_t\in C \mid H_t)=\pi_t(C|H_t)$ for any $C\in\bfrak{B}(\AA)$.
\end{enumerate}
\end{definition}
We have the obvious inclusion $\mathbf{\Pi}\subseteq\tilde{\mathbf{\Pi}}$. 
Moreover, we can associate to each $\pi=(\pi^1,\ldots,\pi^N)\in\mathbf{M}$ the transition kernel (denoted again by~$\pi$) on $\AA$ given $\XX$ defined as
\begin{equation}\label{eq-notation-pidax}
\pi(da|x)=\pi^1(da^1|x)\times\cdots\times\pi^N(da^N|x)\quad\hbox{for $x\in\XX$},
\end{equation}
so that we also have  $\mathbf{M}\subseteq\tilde{\mathbf{M}}$.
For each $\pi\in\tilde{\mathbf{M}}$, we denote by $Q_\pi$ the stochastic kernel on $\XX$ given $\XX$ defined by
\begin{equation*}
Q_{\pi}(D|x)= \int_{\AA} Q(D|x,a)\pi(da|x)
\quad\hbox{for $x\in\XX$ and $D\in\bfrak{X}$}.
\end{equation*}
The compositions of $Q_\pi$ with itself are denoted by $Q^t_{\pi}$ for any $t\ge0$, with the convention that $Q^0_\pi(\cdot|x)$ is the Dirac probability measure concentrated at $x$.

\paragraph{Absorbing games.}

Given a subset of the state space $\Delta\in\bfrak{X}$, we define the hitting time $T_\Delta$ as the measurable function $T_\Delta:\mathbf{H}_\infty\rightarrow\NN\cup\{\infty\}$ given by
$$T_\Delta(x_0,a_0,x_1,a_1,\ldots)=\min\{n\ge0: x_n\in\Delta\},$$
where the $\min$ over the empty set is defined as $+\infty$. Next we propose the some definitions related to the notion of an \textit{absorbing} game.

\begin{definition}
\label{Game-absorbing}
Fix  $\Delta\in\bfrak{X}$ and an initial distribution $\eta\in\bcal{P}(\XX)$.  We say that the game model  
$\mathcal{G}(\eta,\rho)$ is absorbing to $\Delta$  if the conditions (a) and (b) below are satisfied, and we say that 
$\mathcal{G}(\eta,\rho)$ is uniformly absorbing to $\Delta$ if it is absorbing and, in addition, it satisfies condition (c).
\begin{enumerate}
\item[(a).]  For every $(x,a)\in\Delta\times\AA$ we have $Q(\Delta|x,a)=1$ and, besides,  for every  $1\le i\le N$ and $1\le j\le p$ it is $r^i(x,a)=0$ and $c^{i,j}(x,a)=0$;
\item[(b).] For any $\pi\in\tilde{\mathbf{\Pi}}$ the expected hitting time $\mathbb{E}_{\eta,\pi}[T_\Delta]$ is finite.
\item[(c).] We have the following limit: $$\lim_{n\rightarrow\infty} \sup_{\pi\in\tilde{\mathbf{M}}} \sum_{t=n}^\infty\mathbb{P}_{\eta,\pi}\{T_\Delta>t\}=0.$$
\end{enumerate}
\end{definition}

The condition (a) means that, once the state process enters in $\Delta$, it remains in $\Delta$ thereafter at no further reward or cost (related to the functions $r^i$ and $c^i$). The condition (c) can be written in several equivalent ways, for instance:
$$\lim_{n\rightarrow\infty} \sup_{\pi\in\tilde{\mathbf{M}}}  \mathbb{E}_{\eta,\pi}[(T_\Delta-n)^+]=0\quad\hbox{or}\quad
 \sum_{t=0}^n\mathbb{P}_{\eta,\pi}\{T_\Delta>t\}\uparrow \mathbb{E}_{\eta,\pi}[T_\Delta]\ \hbox{uniformly in $\pi\in\tilde{\mathbf{M}}$}.$$ 

Our next result summarizes some important properties. In particular, it is shown that the expected hitting time $\mathbb{E}_{\eta,\pi}[T_\Delta]$ is uniformly bounded over all correlated strategies, which will imply that the set of occupation measures is bounded (see Remark \ref{rem-occupation-measures}(a)), a key element to show the compactness of this set.
\begin{proposition}\label{prop-preliminary} 
Consider a set $\Delta\in\bfrak{X}$ and an initial distribution $\eta\in\bcal{P}(\XX)$.
\begin{itemize} 
\item[(i).] 
If the game model $\mathcal{G}(\eta,\rho)$ is absorbing  to $\Delta$  then
$\sup_{\pi\in\tilde{\mathbf{\Pi}}} \mathbb{E}_{\eta,\pi}[T_\Delta]<\infty$.
\item[(ii).] 
The  family of initial distributions $\eta\in\bcal{P}(\XX)$ for which the game model  
$\mathcal{G}(\eta,\rho)$ is absorbing (respectively, uniformly absorbing) to $\Delta$  is  a convex subset of $\bcal{P}(\XX)$.
\end{itemize}
\end{proposition}
\textbf{Proof.} (i). 
The proof of this item is partially based on some arguments used in \cite[Sections 4.4 and 5.5]{dynkin79} for the special case of a Borel state space.
We will use the following characterization of the set 
$\bcal{S}_\eta=\{\mathbb{P}_{\eta,\pi} \in\bcal{P}(\mathbf{H}_\infty): \pi\in \tilde{\mathbf{\Pi}} \}$ of strategic probability measures. A probability measure
$\mathbb{P}\in\bcal{P}(\mathbf{H}_\infty)$ is in $\bcal{S}_\eta$ if
and only if 
\begin{align}
\mathbb{P}(dx_{0})=\eta(dx_{0})\quad\text{and}\quad\mathbb{P}(dx_0,\ldots,dx_{t},da_{t},dx_{t+1}) =
{\mathbb{P}}(dx_0,\ldots,dx_{t},da_{t}) Q(dx_{t+1}|x_{t},a_{t})
\label{carac-strategic-meas}
\end{align}
for $t\in\NN$, where the above differential notation refers to the marginal of $\mathbb{P}$ on  the corresponding variables.
Let us show the claim by contradiction. So, assume that there exists a sequence $\{\pi_{k}\}_{k\in\NN}$ in $\tilde{\mathbf{\Pi}}$ satisfying 
$\mathbb{E}_{\eta,\pi_{k}}[T_\Delta]\geq 2^k$ for any $k\in\NN$. Consider $\mathbb{P}\in\bcal{P}(\mathbf{H}_\infty)$ defined as
$$\mathbb{P}=\sum_{k\in\NN} \frac{1}{2^{k+1}} \mathbb{P}_{\eta,\pi_{k}}.$$ It is easily seen that $\mathbb{P}$  satisfies \eqref{carac-strategic-meas}. Therefore $\mathbb{P}\in\bcal{S}_\eta$, so that  there exists $\pi\in \tilde{\mathbf{\Pi}}$ with
$\mathbb{P}=\mathbb{P}_{\eta,\pi}$. We have, however, 
$$\mathbb{E}_{\eta,\pi}[T_\Delta]=\int_{\mathbf{H}_\infty} T_\Delta d\mathbb{P}=
\sum_{k\in\NN}  \frac{1}{2^{k+1}}\int_{\mathbf{H}_\infty} T_\Delta d\mathbb{P}_{\eta,\pi_k}=
\sum_{k\in\NN}  \frac{1}{2^{k+1}} \mathbb{E}_{\eta,\pi_{k}}[T_\Delta]=\infty,$$ leading to a contradiction with the condition (b) in Definition \ref{Game-absorbing}.
\par\noindent
(ii). This result, for both the absorbing and the uniformly absorbing cases, is a direct consequence of the fact that $\alpha \mathbb{P}_{\eta,\pi} +(1-\alpha) \mathbb{P}_{\eta',\pi}= \mathbb{P}_{\alpha\eta+(1-\alpha)\eta',\pi}$
for any $\eta,\eta'$ in $\bcal{P}(\mathbf{X})$ and $\pi\in \tilde{\mathbf{\Pi}}$. 
\hfill$\Box$

\paragraph{Equilibria of the game model.}
Given a strategy profile $\pi\in\mathbf{\Pi}$, the total expected payoff of player $i\in\{1,\ldots,N\}$ is
\begin{eqnarray}
R^i(\eta,\pi)=\mathbb{E}_{\eta,\pi} \Big[ \sum_{t=0}^\infty  r^i(X_t,A_t)\Big]=
\mathbb{E}_{\eta,\pi} \Big[ \sum_{0\le t<T_\Delta} r^i(X_t,A_t)\Big]
\in\RR,
\label{expected-payoff}
\end{eqnarray}
and the corresponding total expected cost (for the constraints) is 
\begin{eqnarray}
C^i(\eta,\pi)=\mathbb{E}_{\eta,\pi} \Big[ \sum_{t=0}^\infty  c^i(X_t,A_t)\Big]=
\mathbb{E}_{\eta,\pi} \Big[  \sum_{0\le t<T_\Delta}  c^i(X_t,A_t)\Big]
\in\RR^p.
\label{expected-constraint}
\end{eqnarray}
In Section \ref{sec-assumptions}, we will make assumptions ensuring that $R^i(\eta,\pi)$ and $C^i(\eta,\pi)$ are 
finite ---see  Remark \ref{rem-after-assumption-A}(a)--- for any $\pi\in\mathbf{\Pi}$. For the remainder of this section, we will assume that this is the case.
We say that strategy profile $\pi\in\mathbf{\Pi}$  satisfies the constraint of player $i$  when $C^i(\eta,\pi)\ge\rho_i$.

We propose now the definitions of constrained and unconstrained Nash equilibria.
\begin{definition}  We say that the strategy profile $\pi_*\in\mathbf{\Pi}$ is:
\begin{itemize}
\item[(i).] an unconstrained Nash equilibrium in the class of all strategy profiles if for every $1\le i\le N$ and $\sigma\in\mathbf{\Pi}^i$ we have 
$$  R^i(\eta,(\pi_*^{-i},\sigma))\le R^i(\eta,\pi_*);$$
\item[(ii).] a constrained Nash equilibrium in the class of all strategy profiles if for every $1\le i\le N$ we have $C^i(\eta,\pi_*)\ge\rho^i$ and, in addition, 
$$\forall \sigma\in\mathbf{\Pi}^i,\ C^i(\eta,(\pi_*^{-i},\sigma))\ge\rho^i\ \Rightarrow\ R^i(\eta,(\pi_*^{-i},\sigma))\le R^i(\eta,\pi_*).$$
\end{itemize}
\end{definition}

Next we introduce the standard Slater condition.  It states that whatever policies the other players use, player $i$ can find a policy so as to satisfy his own constraints.
\begin{definition}\label{def-slater}
We say that the game model $\mathcal{G}(\eta,\rho)$ satisfies the Slater condition when, for each strategy profile $\pi\in\mathbf{\Pi}$ and any player $1\le i\le N$, there exists $\sigma^i\in\mathbf{\Pi}^i$ such that $C^i(\eta,(\pi^{-i},\sigma^i))>\rho^i$.  
\end{definition}

\subsection{Assumptions and Young measures}
\label{sec-assumptions}
In this section we will formulate the assumptions we will need in the sequel. We will also introduce the notion of Young measure.
{We present three sets of basic assumptions that slightly differ from each other. It is easy to notice that assumption \ref{Assump} is weaker than
\ref{Assump-2}, which is itself weaker than \ref{Assump-3}.
In what follows, many results will be proved for an arbitrary initial distribution in $\bcal{P}_\lambda(\XX)$ which will be noted generically by $\eta$.} 
\begin{hypot}
\label{Assump}
\item \mbox{ } 
Consider the game model $\mathcal{G}(\eta,\rho)$  with initial distribution $\eta\in\bcal{P}(\XX)$ and constraint constants $\rho\in\RR^{pN}$.
{We say that $\mathcal{G}(\eta,\rho)$ satisfies Assumption \ref{Assump} when there exist a probability measure $\lambda\in\bcal{P}(\XX)$
and a set $\Delta\in\bfrak{X}$ for which the following conditions hold:}
\begin{assump}
\item\label{Assump-absorbing}  The game model $\mathcal{G}(\eta,\rho)$ is absorbing to~$\Delta$.
\item \label{Assump-Slater} The game model $\mathcal{G}(\eta,\rho)$ satisfies the Slater condition.
\item \label{Assumption-state} The $\sigma$-algebra $\bfrak{X}$ is countably generated.
\item \label{Assumption-control} For each player $i\in\{1,\ldots,N\}$, the action set $\AA^i$  is compact 
and the correspondence from $\XX$ to $\AA^i$ given by $x\mapsto \AA^i(x)$ is weakly measurable with nonempty compact values.
 \item \label{Assumption-reward} For each player $i\in\{1,\ldots,N\}$, we have that $r^i$ and $c^i$  are bounded Carath\'eodory functions, that is, $r^i\in\car_b(\XX\times\AA,\RR)$ and  $c^i\in\car_b(\XX\times\AA,\RR^p)$. Let $\mathbf{r}>0$ be a componentwise bound  for all the $r^i$ and $c^i$.
\item \label{Assumption-transition-Q} There exists a measurable density  function $q: \XX\times\XX\times\AA\rightarrow\RR^+$  such that  for each $B\in\bfrak{X}$  and  $(x,a)\in\XX\times\AA$ we have
$$
Q(B|x,a)=\int_B  q(y,x,a)\lambda(dy)
\quad\hbox{and}\quad
 \lim_{n\rightarrow\infty} \int_{\XX} |q(y,x,a_n)-q(y,x,a)| \lambda(dy)=0
\label{Assumption-transition-Q-Nowak}
$$
whenever $a_n\rightarrow a$ in $\AA$.
\item\label{Assump-absorbing-lambda}  The game model $\mathcal{G}(\lambda,\rho)$ is absorbing to~$\Delta$.
\item \label{Assump-Abs-Cont}  The initial distribution $\eta$ satisfies $\eta\ll\lambda$.
\end{assump}
\end{hypot}

\begin{hypotprim}
\label{Assump-2}
\item \mbox{ } 
Consider the game model $\mathcal{G}(\eta,\rho)$  with initial distribution $\eta\in\bcal{P}(\XX)$ and constraint constants $\rho\in\RR^{pN}$.
We say that $\mathcal{G}(\eta,\rho)$ satisfies Assumptions \ref{Assump-2} when it satisfies Assumption \ref{Assump}  except for \ref{Assump-absorbing} which is replaced by the following stronger condition:
\begin{assumpprim}
\item\label{Assump-absorbing-uniformly}  The game model $\mathcal{G}(\eta,\rho)$ is uniformly absorbing to~$\Delta$.
\end{assumpprim}
\end{hypotprim}

\begin{hypotcal}
\label{Assump-3}
\item \mbox{ } 
Consider the game model $\mathcal{G}(\eta,\rho)$  with initial distribution $\eta\in\bcal{P}(\XX)$ and constraint constants $\rho\in\RR^{pN}$.
We say that $\mathcal{G}(\eta,\rho)$ satisfies Assumptions \ref{Assump-3} when it satisfies Assumption \ref{Assump-2}  except for \ref{Assump-absorbing-lambda} which is replaced by the following stronger condition:
\begin{assumpcal}
\setcounter{assumpcal}{6}
\item\label{Assump-absorbing-lambda-uniformly}  The game model $\mathcal{G}(\lambda,\rho)$ is uniformly absorbing to~$\Delta$.
\end{assumpcal}
\end{hypotcal}
The Assumptions \ref{Assump}, \ref{Assump-2} and \ref{Assump-3} will be discussed in the next section.

\paragraph{The space of Young measures $\bcal{Y}$.} 
Now  we introduce the notion of Young measure in order to endow the spaces $\mathbf{M}$ and $\tilde{\mathbf{M}}$ of stationary Markov strategies with a suitable metric. 
Note also that the last assumption we will need in this paper will be formulated at the end of this section, in terms of continuity properties of functions defined on these Young measure spaces.
To do this we will rely on the reference probability measure $\lambda$ on the state space $\mathbf{X}$ introduced in Assumption~\ref{Assumption-transition-Q}. 

Recall that, for  the noncooperative game model, the class of stationary Markov profiles of the players is  $\mathbf{M}^1\times\ldots \times\mathbf{M}^N$.
Given a player $1\le i\le N$, we consider in $\mathbf{M}^{i}$ the following equivalence relation: for  $\phi,\varphi\in\mathbf{M}^i$ we say that 
$$\phi\sim\varphi\quad\hbox{when}\quad
\phi(\cdot|x)=\varphi(\cdot|x)\quad\hbox{for $\lambda$-almost every $x\in\XX$}.$$  We will denote by  $\bcal{Y}^i$ the corresponding family of equivalence classes,
which will be referred to as \textit{Young measures}. The set $\bcal{Y}^i$ of Young measures is equipped with the narrow (stable) topology: it is the coarsest topology that makes the mappings 
$$\pi^i\mapsto \int_{\XX}\int_{\AA^i}f(x,a^i) \pi^i (da^i|x)\lambda(dx)$$ continuous for any $f$ which is a Carath\'eodory function on $\XX\times\AA^i$ bounded in $L^1$; more precisely, this means that  $f\in\car(\XX\times\AA^i,\RR)$ is such that for some $F\in L^1(\XX,\bfrak{X},\lambda)$
we have $|f(x,a^i)| \leq F(x)$ for every $(x,a^{i})\in\XX\times\AA^i$; see, e.g.,  \cite[Theorem 2.2]{balder88}.
Using \cite[Lemma 1]{balder91}, it follows that $\bcal{Y}^i$ is a compact metric space for the narrow topology. 
We also define
 $$\bcal{Y}=\bcal{Y}^1\times\ldots\times \bcal{Y}^N,$$ which is endowed with the product topology, and it is therefore a compact metric space as well. 
 
\paragraph{Young measures $\bcal{Y}$  and Markov strategies $\mathbf{M}$.}
We will say that two Markov strategies of the noncooperative game model $\pi=(\phi^1,\ldots,\phi^N)$ and $\pi'=(\varphi^1,\ldots,\varphi^N)$ in $\mathbf{M}$ are in the same equivalence class of Young measures whenever  $\phi^i\sim\varphi^i$ for every $1\le i\le N$, and we will write $\pi\sim\pi'$ as well. In this case, since the initial distribution $\eta$ is absolutely continuous with respect to $\lambda$, and since the transition of the system has a density with respect to $\lambda$, it is easily seen \cite[Lemma 2.2]{Dufour-Prieto23} that both strategies yield the same strategic probability measure, that is, $\mathbb{P}_{\eta,\pi}=\mathbb{P}_{\eta,\pi'}$. Therefore, $\pi$ and $\pi'$ are indistinguishable since they are driven by the same strategic probability measures and, besides, they also have the same  costs and rewards since those are defined from the corresponding strategic measures. Hence, in the sequel \textit{we shall identify the set of Young measures~$\bcal{Y}$ with the family of Markov profiles $\mathbf{M}$ of the players}. 

\paragraph{Consistence of the notation.}
Given $\pi\in\mathbf{M}$
and a function $f\in\car_b(\XX\times\AA)$,   define the measurable function $f_\pi$ on $\XX$ as
$$f_\pi(x)=\int_\AA f(x,a)\pi(da|x)\quad\hbox{for $x\in\XX$}.$$
If $\pi'\in\mathbf{M}$ is such that $\pi'\sim\pi$ then $f_{\pi'}=f_\pi$ with $\lambda$-probability one, and so they both belong to the same equivalence class in $L^\infty(\XX,\bfrak{X},\lambda)$. Therefore, it is consistent to define the function $f_\pi\in L^\infty(\XX,\bfrak{X},\lambda)$ for $\pi\in\bcal{Y}$.

Suppose that 
 $v\in L^\infty(\XX,\bfrak{X},\lambda)$ and $\pi\in\mathbf{M}$. We have
$$Q_\pi v(x)=\int_\AA \int_\XX v(y)q(y,x,a)\lambda(dy)\pi(da|x)\quad\hbox{for $x\in\XX$}.$$
Hence, the above integral does not depend on the representative $v$ chosen in $L^\infty(\XX,\bfrak{X},\lambda)$; and, in addition, if $\pi'\in\mathbf{M}$ is in the same equivalence class of Young measures as $\pi$, then $Q_{\pi}v=Q_{\pi'}v$ with $\lambda$-probability one, and so both $\pi$ and $\pi'$ yield the same element in $L^\infty(\XX,\bfrak{X},\lambda)$. Consequently, the notation 
$Q_\pi v\in L^\infty(\XX,\bfrak{X},\lambda)$ for $v\in L^\infty(\XX,\bfrak{X},\lambda)$  and  $\pi\in\bcal{Y}$ is consistent.
The same applies for the successive compositions $Q^t_\pi v$ of the stochastic kernels for $t\ge0$.
Note also  that Assumption~\ref{Assumption-transition-Q}  implies, in particular, that $Qv$ is well defined, meaning that  $Qv=Qv'$ whenever $v$ and $v'$ belong to the same equivalence class in  $L^{\infty}(\XX,\bfrak{X},\lambda)$.

\paragraph{Young measures $\tilde{\bcal{Y}}$  and Markov strategies $\tilde{\mathbf{M}}$.}
We recall that a  stationary correlated Markov strategy $\pi\in\tilde{\mathbf{M}}$ is given by a stochastic kernel $\pi$ on $\AA$ given $\XX$ with $\pi(\AA(x)|x)=1$ for every $x\in\XX$.  As before, we can identify kernels $\pi,\pi'\in\tilde{\mathbf{M}}$ such that $\pi(\cdot|x)=\pi'(\cdot|x)$  $\lambda$-a.s. on $\XX$ (written $\pi\sim\pi'$) and then define the set $\tilde{\bcal{Y}}$ of Young measures as the corresponding equivalence classes. The  associated narrow topology is the coarsest one that makes the mappings
$\pi\mapsto \int_{\XX}\int_{\AA} f(x,a)\pi(da|x)\lambda(dx)$ continuous for every $f\in\car(\XX\times\AA,\RR)$ bounded by a $\lambda$-integrable function. Again, 
 $\tilde{\bcal{Y}}$ is a compact metric space with its narrow topology.

 Observe that the equivalence relation of $\bcal{Y}$ is compatible with that of $\tilde{\bcal{Y}}$ meaning that, given $\pi,\pi'\in\mathbf{M}\subseteq\tilde{\mathbf{M}}$,  we have $\pi\sim\pi'$ in $\mathbf{M}$ if and only if $\pi\sim\pi'$ in $\tilde{\mathbf{M}}$.
 Notations such as $f_\pi$  and $Q_\pi^t v$ for $\pi\in\tilde{\bcal{Y}}$, $f\in\car_b(\XX\times\AA)$, $v\in L^\infty(\XX,\bfrak{X},\lambda)$, and $t\ge0$ are consistent as well. 
Invoking the same previous arguments, we shall hereafter identify the space $\tilde{\bcal{Y}}$ of Young measures and the class $\tilde{\mathbf{M}}$ of correlated Markov strategies of the players.

\noindent
We make the following very important remark.
\begin{remark}
\label{rem-trace}
Although it is true that $$\bcal{Y}^1\times\ldots\times \bcal{Y}^N=\bcal{Y}\subseteq\tilde{\bcal{Y}},$$  it turns out that the trace of the narrow topology of $\tilde{\bcal{Y}}$ on $\bcal{Y}$ does not coincide, in general, with the product topology of the $\bcal{Y}^i$. Namely, suppose that $\{\pi_n\}$ and $\pi$ are in $\bcal{Y}$. If $\pi_n\rightarrow\pi$ in $\tilde{\bcal{Y}}$ then it is easy to verify that $\pi_n\rightarrow\pi$ in $\bcal{Y}$. The converse: $\pi_n\rightarrow\pi$ in ${\bcal{Y}}$ implies $\pi_n\rightarrow\pi$ in $\tilde{\bcal{Y}}$, however,
is not necessarily true.
Therefore, to fix the terminology, by convergence in $\bcal{Y}$ we shall refer to convergence in the product topology of the $\bcal{Y}^i$, whereas convergence in $\tilde{\bcal{Y}}$ will mean convergence in the narrow topology of $\tilde{\bcal{Y}}$. 
\end{remark}

We now introduce an additional condition that will allow us to establish continuity results for the game model. 
These continuity conditions will be expressed in terms of functions defined on the set of Markov strategies $\bcal{Y}$.
In the next section, we will propose sufficient conditions for Assumption~\ref{Assumption-B} below.

\begin{hypot}\label{Assumption-B}
\item{} 
We say that the game model $\mathcal{G}(\eta,\rho)$ satisfies Assumption \ref{Assumption-B} when
the following mappings, defined on $\bcal{Y}$ and taking  values in $L^{\infty}(\XX,\bfrak{X},\lambda)$,
$$\pi\mapsto r^i_\pi\, ,\quad \pi\mapsto c^{i,j}_\pi\, ,\quad\hbox{and}\quad \pi\mapsto Q_\pi v$$
are continuous for any $1\le i\le N$, $1\le j\le p$, and $v\in L^{\infty}(\XX,\bfrak{X},\lambda)$.
\end{hypot}

\section{Main results}\label{sec-3}
In this section we will present our two main results. The first one shows the existence of a Markovian equilibrium in the case where $\lambda$ is absolutely continuous with respect to the initial distribution $\eta$ and under assumptions \ref{Assump-2} and \ref{Assumption-B}.  The second result relaxes the condition $\lambda\ll \eta$ but requires strengthening Assumption \ref{Assump-2} to \ref{Assump-3}.

\begin{proposition} 
\label{th-main-0}
Suppose that we are given an initial distribution $\eta\in\bcal{P}(\XX)$ and constraint constants $\rho\in\RR^{Np}$ such that $\lambda\ll \eta$, and that the game $\mathcal{G}(\eta,\rho)$ satisfies Assumptions~\ref{Assump-2} and \ref{Assumption-B}. 
Both the constrained and the unconstrained game have a stationary Markov profile which is a Nash equilibrium in the class of all strategy profiles.
\end{proposition}
\textbf{Proof.}
See Section \ref{Proof-Main1}
$\hfill$ $\Box$

\begin{theorem} 
\label{th-main}
Suppose that we are given an initial distribution $\eta\in\bcal{P}(\XX)$ and constraint constants $\rho\in\RR^{Np}$ such that the game $\mathcal{G}(\eta,\rho)$ satisfies
Assumptions~\ref{Assump-3} and \ref{Assumption-B}. 
Both the constrained and the unconstrained game have a stationary Markov profile which is a Nash equilibrium in the class of all strategy profiles.
\end{theorem}
\textbf{Proof.}
See Section \ref{Proof-Main2}
$\hfill$ $\Box$

\begin{remark}
\label{elements-of-proof}
The proof of Theorem \ref{th-main} will proceed in several steps.
First, we need to consider  the case where the reference probability measure $\lambda$ is absolutely continuous with respect to the initial distribution $\eta$ and show the existence of a Markovian equilibrium under assumptions \ref{Assump-2} and \ref{Assumption-B}, which is precisely Proposition \ref{th-main-0}.
Then, in a second step, we will drop the hypothesis $\lambda\ll\eta$ by considering a sequence of game models $\mathcal{G}(\eta_{n},\rho_{n})$
that satisfies Assumption \ref{Assump-2} and \ref{Assumption-B} for  $\eta_{n}$ given by a convex combinations of $\eta$ and $\lambda$:
\begin{equation}\label{eq-def-zeta-n}
\eta_n=\frac{n}{n+1}\eta+\frac{1}{n+1}\lambda
\end{equation}
and for suitably defined constraint constants $\rho_{n}$.
On the basis of Proposition \ref{th-main-0}, this yields the existence of a constrained Nash equilibrium $\hat{\pi}_{n}\in\bcal{Y}$ for the game model $\mathcal{G}(\eta_n,\rho_n)$.
Finally, we will prove Theorem  \ref{th-main} by showing that a limit point of the sequence $\{\hat{\pi}_{n}\}\in\bcal{Y}$ provides a Markovian equilibrium for the game model
$\mathcal{G}(\eta,\rho)$.
\end{remark}

Some comments regarding the assumptions are in order now.
\begin{remark}
\label{comment-assumption-A}
\begin{enumerate}
\item[(a).]
In the context of an absorbing game model $\mathcal{G}(\eta,\rho)$, the assumptions \ref{Assump-absorbing}-\ref{Assumption-transition-Q} and \ref{Assump-Abs-Cont} are conditions classically met in the literature, see for example \cite{dufour22} for the special case of a discounted model.
Assumption \ref{Assump-absorbing-lambda} is a key technical condition which will allow to show very important properties of the absorbing model
by showing in particular that a measure in $\bcal{M}^+(\XX\times\AA)$  is an occupation measure if and only if it satisfies the characteristic equations (see 
items (i) and (iv) of Proposition \ref{prop-linear-equation-occupation}). 
From this point of view, all the conditions of Assumption \ref{Assump} are very natural.
\item[(b).] The proof of the existence of a Markovian noncooperative equilibrium relies on the fact that the set of the marginals on $\XX\times\AA^i$
of the occupation measure is a compact space in order to use the Kakutani-Fan-Glicksberg fixed point Theorem. It will be shown in Proposition \ref{Compactness-set-D} that, under Assumption \ref{Assump}, this set is compact if and only if  the model is uniformly absorbing.
This is why it is necessary to replace Assumption \ref{Assump} by Assumption \ref{Assump-2} to prove Proposition~\ref{th-main-0}.
Now regarding the proof of Theorem \ref{th-main}, a key step is to show $\mathcal{G}(\eta_n,\rho_n)$ is uniformly absorbing as explained in Remark \ref{elements-of-proof}.
It is, therefore, necessary to reinforce the absorbing condition of $\mathcal{G}(\lambda,\rho)$ by assuming that $\mathcal{G}(\lambda,\rho)$ is uniformly absorbing, which leads to replace Assumption \ref{Assump-2} by Assumption \ref{Assump-3} in the statement of the theorem \ref{th-main}.
\end{enumerate}
\end{remark}

Let us now describe some consequences of these assumptions.
\begin{remark}\label{rem-after-assumption-A}
\begin{enumerate}
\item[(a).] {The functions $r^i$ and $c^i$ being bounded under \ref{Assumption-reward} and the game model 
$\mathcal{G}(\eta,\rho)$ being absorbing to~$\Delta$ under \ref{Assump-absorbing} (and also \textit{a fortiori} under \ref{Assump-absorbing-uniformly}) we have that
$R^i(\eta,\pi) \leq\mathbf{r} \sup_{\pi\in\Pi}\EE_{\eta,\pi}[T_\Delta]$ and $C^i(\eta,\pi)\leq\mathbf{r} \sup_{\pi\in\Pi}\EE_{\eta,\pi}[T_\Delta] \1$ (see equations \eqref{expected-payoff} and \eqref{expected-constraint}), which  are therefore finite for any $\pi\in\Pi$ and each player $i$ by using Proposition \ref{prop-preliminary}.}
\item[(b).] The Slater condition \ref{Assump-Slater} implies  that $\lambda(\Delta^c)>0$. Indeed, otherwise  we would have $\eta(\Delta)=1$ and so the process would always remain with probability one in $\Delta$ and the corresponding reward functions would be all zero. 
In this very particular case where $\lambda(\Delta^c)=0$, the problem is degenerate and all Markov strategies are noncooperative equilibria.
\item[(c).] By Assumption \ref{Assumption-control}, the correspondences $x\mapsto \AA^i(x)$  are measurable \cite[Lemma 18.2]{aliprantis06} and they have measurable graph  \cite[Theorem 18.6]{aliprantis06}. Therefore,  the following sets are measurable:
\begin{align*}
&\KK^{i} =\{ (x,a^{i})\in\XX\times\AA^{i} : a^{i}\in\AA^{i}(x) \}\in \bfrak{X}\otimes\bfrak{B}(\AA^{i}) \quad\hbox{for $1\le i\le N$, and} 
\\
 &\KK=\{ (x,a)\in\XX\times\AA : a\in\AA(x)\}\in \bfrak{X}\otimes\bfrak{B}(\AA).
\end{align*}
The  Kuratowski-Ryll-Nardzewski selection theorem \cite[Theorem 18.13]{aliprantis06} yields the existence of measurable selectors for $x\mapsto \AA_i(x)$ for each $1\le i\le N$. In particular, $\mathbf{M}^i$ is nonempty, and so are all the sets of strategies defined so far: $\mathbf{\Pi}^i$, $\mathbf{\Pi}$,  $\tilde{\mathbf{\Pi}}$, $\mathbf{M}$, and $\tilde{\mathbf{M}}$.
\end{enumerate}
\end{remark}

We conclude this section by proposing sufficient conditions yielding the continuity properties stated in Assumption~\ref{Assumption-B}. 
We show that some game models classically studied in the literature satisfy our hypotheses, such as countable state space models or ARAT-type models.
It is also important to note that the expected discounted models are a special case of the total expected absorbing models. In this way, our results here generalize those in \cite{dufour22,Jaskiewicz-Nowak-AMO-22}.
\begin{corollary}
\label{Cond-sufficient-Assumption-B}
Suppose that the game model $\mathcal{G}(\eta,\rho)$  satisfies Assumption \ref{Assump}. Under any of the conditions (i) and (ii) below, Assumption \ref{Assumption-B} is satisfied.
\begin{itemize}
\item[(i).] The state space $\XX$ is countable.
\item[(ii).] The game model has ARAT structure (additive reward, additive transition) meaning that:
\begin{itemize}
\item[(a).] (Additive reward.) For any $1\le i\le N$, $1\le j\le p$, and $1\le l\le N$, there exist functions $r^i_l$ and $c^{i,j}_l$ in $\car_b(\XX\times\AA^l)$ such that 
$$r^i(x,a^1,\ldots,a^N)=\sum_{l=1}^N r^i_l(x,a^l)\quad\hbox{and}\quad
c^{i,j}(x,a^1,\ldots,a^N)=\sum_{l=1}^N c^{i,j}_l(x,a^l)$$
for any $(x,a)\in\XX\times\AA$.
\item[(b).] (Additive transition.) There exist nonnegative measurable functions $q^l:\XX\times\XX\times\AA^l\rightarrow\RR$ such that 
$$Q(B|x,a^1,\ldots,a^N)=\sum_{l=1}^N \int_B q^l(y,x,a^l)\lambda(dy)\quad\hbox{for $B\in\bfrak{X}$ and $(x,a^1,\ldots,a^N)\in\XX\times\AA$}$$
with, in addition, $\lim_{n\rightarrow\infty}\int_\XX | q^l(y,x,a^l_n)-q^l(y,x,a)|\lambda(dy)=0$ for any $x\in\XX$ whenever $a_n^l\rightarrow a^l$ as $n\rightarrow\infty$ in $\AA^l$.
\end{itemize}
\end{itemize}
\end{corollary}
\textbf{Proof.} See Section \ref{sec53}.\hfill$\Box$

\begin{remark}
The absolute continuity condition in Assumption \ref{Assump-Abs-Cont} is not restrictive with respect to Assumptions \ref{Assump}, \ref{Assump-2} and \ref{Assump-3}. 
Indeed, if it were not true that $\eta\ll\lambda$, then we would  consider  the reference probability measure $\bar{\lambda}=(\eta+\lambda)/2$.
Clearly we have $\eta\ll\bar\lambda$, while   
 the function $\bar{q}(y,x,a)=q(y,x,a)(d\lambda/d\bar{\lambda})(y)$  would satisfy Assumption \ref{Assumption-transition-Q}.  
Finally, regarding Assumptions \ref{Assump-absorbing-lambda} or \ref{Assump-absorbing-lambda-uniformly}, the convexity property in Proposition \ref{prop-preliminary}(ii)
ensures that the game model $\mathcal{G}(\bar\lambda,\rho)$ is (respectively, uniformly) absorbing to~$\Delta$
if $\mathcal{G}(\lambda,\rho)$ and $\mathcal{G}(\eta,\rho)$ are (respectively, uniformly) absorbing.
However, changing the measure $\lambda$ to $\bar{\lambda}$ may affect Assumption \ref{Assumption-B}. 
Nevertheless, it is important to emphasize that conditions (i) and (ii) of Corollary \ref{Cond-sufficient-Assumption-B} implying Assumption \ref{Assumption-B} are not affected by a change of the reference probability measure $\lambda$.
\end{remark}

\paragraph{A note on discounted games.}
As mentioned in \cite[p. 132]{Feinberg-Rothblum}, a $\beta$-discounted model can be transformed into an equivalent absorbing model just by adding an isolated absorbing cemetery state $x_\Delta$ with a single available action $a_\Delta$ at no reward or cost. In this way,  the new state space is $\XX'=\XX\cup \{x_\Delta\}$ and    the transitions of the system are given by
$$Q'(B|x,a)=\begin{cases} \beta Q(B|x,a)\ \hbox{when $B\subseteq\XX$}\\ 1-\beta\ \hbox{if $B=\{x_\Delta\}$}
\end{cases}$$
for $(x,a)\in\XX\times\AA$ and $Q'(\{x_\Delta\}|x_\Delta,a_\Delta)=1$.
The reference probability measure would be
$$\lambda'(B)=\beta\lambda(B\cap\XX)+(1-\beta)\delta_{\{x_\Delta\}}(B)\quad\hbox{for measurable $B\subseteq\XX'$}.$$
It is then easily seen that this game model is uniformly absorbing to $\{x_\Delta\}$ for any initial distribution.
As a consequence, Assumptions \ref{Assump-absorbing}, \ref{Assump-absorbing-lambda}, \ref{Assump-absorbing-uniformly}
and \ref{Assump-absorbing-lambda-uniformly}
can be dropped, which makes Assumptions \ref{Assump}, \ref{Assump-2} and \ref{Assump-3} equivalent. 
This shows that under Assumptions \ref{Assump} and \ref{Assumption-B}, we obtain the existence of constrained and unconstrained equilibria for discounted games.
As a consequence, using Corollary \ref{Cond-sufficient-Assumption-B} we would obtain the results in \cite{dufour22} and \cite{Jaskiewicz-Nowak-AMO-22}.

\section{Occupation measures} 
\label{section-Assumptions-YoungMeasures}
\subsection{Occupation measures and their topological properties}
\textit{Throughout this subsection we shall assume that we are given a game model $\mathcal{G}(\eta,\rho)$ with initial distribution $\eta\in\bcal{P}_\lambda(\XX)$ and constraint constant $\rho\in\RR^{pN}$ that satisfies Assumption \ref{Assump}.}

First of all we state some useful properties of the kernel $\mathbb{I}_{\Delta^c}$ on $\XX$ given~$\XX$ which was defined in Section \ref{sec-1}.
\begin{lemma}\label{lemma-Identity-absorbing} 
The kernel $\mathbb{I}_{\Delta^c}$ on $\XX$  given $\XX$ satisfies the following properties.
Given any $\pi\in\tilde{\bcal{Y}}$, $f\in L^\infty(\XX,\bfrak{X},\lambda)$, and
$\mu\in\bcal{M}^+(\XX)$:
\begin{itemize}  
\item[(i).] $Q_\pi\mathbb{I}_{\Delta^c}(B|x)=Q_\pi(B\cap\Delta^c|x)$ for $B\in\bfrak{X}$ and $x\in\XX$.
\item[(ii).] $\mathbb{I}_{\Delta^c}Q_\pi(B|x)=Q_\pi(B|x)\mathbf{I}_{\Delta^c}(x)$ for $B\in\bfrak{X}$ and $x\in\XX$.
\item[(iii).] $\mathbb{I}_{\Delta^c}Q_\pi\mathbb{I}_{\Delta^c}=Q_\pi\mathbb{I}_{\Delta^c}$ and, as a consequence, $\mathbb{I}_{\Delta^c}(Q_\pi\mathbb{I}_{\Delta^c})^t=Q_\pi^t\mathbb{I}_{\Delta^c}=(Q_\pi\mathbb{I}_{\Delta^c})^t$ for any $t\ge1$.
\item[(iv).] $(\mathbb{I}_{\Delta^c}f)(x)=f(x)\mathbf{I}_{\Delta^c}(x)$ for $x\in\XX$.
\item[(v).] $\mu\mathbb{I}_{\Delta^c}(B)=\mu(B\cap\Delta^c)$, which can be also written $\mu\mathbb{I}_{\Delta^c}=\mu_{\Delta^c}$.
\end{itemize}
\end{lemma}

Next we propose the definition of the occupation measure induced by a correlated strategy of the players.
This definition can be specialized to noncooperative strategy profiles.
We recall that we are making the convention that the sum over an empty set is zero.

\begin{definition}\label{def-occupation-measure}
Given any strategy profile $\pi\in\tilde{\mathbf{\Pi}}$, 
the  \textit{occupation measure}   $\mu_{\eta,\pi}\in\bcal{M}^{+}(\XX\times\AA)$ for the initial distribution $\eta\in\bcal{P}_\lambda(\XX)$ is defined, for measurable sets  $B\in\bfrak{X}$ and $D^i\in\mathfrak{B}(\AA^i)$ for $1\le i\le N$, as
\begin{eqnarray*}
\mu_{\eta,\pi}(B\times D^1\times\ldots\times D^N)
&=& \mathbb{E}_{\eta,\pi}\Big[ \sum_{t=0}^\infty \mathbf{I}_{\{T_\Delta>t\}}\cdot\mathbf{I}_{\{X_t\in B, A^1_t\in D^1,\ldots, A^N_t\in D^N\}}\Big].
\end{eqnarray*}
We introduce the notations $
\tilde{\bcal{O}}_\eta=\{\mu_{\eta,\pi}:\pi\in\tilde{\mathbf{\Pi}}\}$ and 
$\bcal{O}_\eta=\{\mu_{\eta,\pi}:\pi\in\mathbf{\Pi}\}$,
with $\bcal{O}_\eta\subseteq \tilde{\bcal{O}}_\eta$.
\end{definition}

Some important comments concerning this definition are in order.
\begin{remark}\label{rem-occupation-measures}
\begin{enumerate}
\item[(a).] 
Note that $\mu_{\eta,\pi}(\XX\times\AA)=\mathbb{E}_{\eta,\pi}[T_\Delta]$ for $\pi\in\tilde{\mathbf{\Pi}}$, which is finite because $\mathcal{G}(\eta,\rho)$ is absorbing to $\Delta$. Moreover, by Proposition \ref{prop-preliminary}, we have
$$\sup_{\pi\in\tilde{\mathbf{\Pi}}} \mu_{\eta,\pi}(\XX\times\AA)<\infty,$$
which is usually referred to as $\tilde{\bcal{O}}_\eta$ being bounded. Also observe that, by construction of the process, $\mu_{\eta,\pi}(\KK^c)=0$.
Clearly, the set $\bcal{O}_\eta$ of occupation measures of the noncooperative game inherits the same properties.
\item[(b).] Observe that, although the process will eventually visit the set $\Delta$ ---it might even be $\eta(\Delta)>0$---  we have $\mu_{\eta,\pi}^\XX(\Delta)=0$. This is because, by its definition, the occupation measure ``does not count'' visits to $\Delta$. In fact,  the $\XX$-marginal of the occupation measure is given by 
\begin{equation}\label{eq-marginal-occupation-measure}
\mu_{\eta,\pi}^\XX(B)=\mathbb{E}_{\eta,\pi}\big[\sum_{t=0}^\infty  \mathbf{I}_{\{T_\Delta>t\}}\cdot\mathbf{I}_{\{X_t\in B\}}\big]=
\sum_{t=0}^\infty \mathbb{P}_{\eta,\pi}\{X_t\in B-\Delta\}\quad \hbox{for $B\in\bfrak{X}$}
\end{equation}
because we have $\{X_t\in B,T_\Delta>t\}=\{X_t\in B-\Delta\}$.
\item[(c).]
 It follows directly from  Definitions \ref{Game-absorbing} and \ref{def-occupation-measure}, and Assumption  \ref{Assumption-reward}  that the total expected payoffs of the strategy profile $\pi\in\mathbf{\Pi}$  for the initial distribution $\eta\in\bcal{P}_\lambda(\XX)$ equal
\begin{equation}\label{eq-R-as-a-sum}
R^i(\eta,\pi)=\int_{\XX\times\AA} r^id\mu_{\eta,\pi}\quad\hbox{and}\quad
C^i(\eta,\pi)=\int_{\XX\times\AA} c^id\mu_{\eta,\pi}\quad\hbox{for $1\le i\le N$}.
\end{equation}
\item[(d).] Regarding Markov strategies in $\mathbf{M}$ or $\tilde{\mathbf{M}}$, since their occupation measure is defined based on the corresponding strategic probability measures, if follows that two Markov strategies in the same equivalence class of  $\bcal{Y}$ or $\tilde{\bcal{Y}}$ yield the same occupation measure. So, the notation $\mu_{\eta,\pi}$ for $\pi \in\bcal{Y}$ or $\pi\in\tilde{\bcal{Y}}$ is consistent.
\end{enumerate}
\end{remark}

Let us first show the following technical results before deriving properties on occupation measures. 
\begin{lemma}
\label{prop-survival-TDelta}
Let $\boldsymbol{\Gamma}$ be an arbitrary subset of $\tilde{\mathbf{M}}$ and let $\{h_{\pi}\}_{\pi\in\boldsymbol{\Gamma}}$ be a family of non-negative functions in
$L^{1}(\XX,\bfrak{X},\lambda)$ which are uniformly $\lambda$-integrable. Under these conditions,
$$
\lim_{t\rightarrow\infty} \sup_{\pi\in \boldsymbol{\Gamma}} \int_\XX Q^t_\pi(\Delta^c|x)h_\pi(x)\lambda(dx)= 0.$$
\end{lemma}
\textbf{Proof.}
Consider a fixed arbitrary $\epsilon>0$. By the uniform integrability hypothesis, there exists $c_{\epsilon}>0$ such that
$$ \sup_{\pi\in \boldsymbol{\Gamma}} \int_{\{x\in\XX: h_{\pi}(x)>c_{\epsilon}\}} h_{\pi}(x) \lambda (dx)\leq \epsilon.$$
Therefore, for any  $\pi\in\mathbf{\Gamma}$  and $t\ge1$
$$
 \int_\XX Q^t_\pi(\Delta^c|x)h_\pi(x)\lambda(dx)
 \leq \epsilon+c_{\epsilon} \mathbb{P}_{\lambda,\pi}\{T_{\Delta}>t\}
\leq \epsilon+ \frac{c_{\epsilon}}{t} \cdot\mathbb{E}_{\lambda,\pi}[T_{\Delta}].
$$
From Assumption \ref{Assump-absorbing-lambda} and applying Proposition \ref{prop-preliminary}(i) we have that $\sup_{\pi\in\boldsymbol{\Gamma}} \mathbb{E}_{\lambda,\pi}[T_\Delta]<\infty$. Hence, we choose $t$ large enough so as to obtain that 
$\sup_{\pi\in \boldsymbol{\Gamma}} \int_\XX Q^t_\pi(\Delta^c|x)h_\pi(x)\lambda(dx)<2\epsilon$, 
and the result follows.
$\hfill$ $\Box$

\begin{lemma}
\label{prop-linear-equation-occupation-unicity}
Given any $\pi\in \tilde{\mathbf{M}}$, the measure $\gamma \in \bcal{M}^+(\XX)$ defined as
\begin{align}
\label{expression-marginal}
\gamma=
\eta \sum_{k=0}^{\infty} Q_{\pi}^{k}\mathbb{I}_{\Delta^c}
\end{align}
satisfies $\gamma\ll\lambda$  and it is the unique solution  of  the equation
\begin{align}
\label{eq-marginal}
\xi=(\eta +  \xi Q_{\pi})\mathbb{I}_{\Delta^c} \quad\hbox{for $\xi \in \bcal{M}^+(\XX)$}.
\end{align}
Moreover, $\gamma=\mu_{\eta,\pi}^\XX$.
\end{lemma}
\textbf{Proof.} First of all, observe that $\gamma$ defined in \eqref{expression-marginal} is indeed in $\bcal{M}^+(\XX)$ because
$$\gamma(\XX)=\sum_{k=0}^\infty ( \eta Q^k_\pi \mathbb{I}_{\Delta^c})(\XX)=\sum_{k=0}^\infty (\eta Q^k_\pi)(\Delta^c)=\sum_{k=0}^\infty \mathbb{P}_{\eta,\pi}\{T_\Delta>k\}=\mathbb{E}_{\eta,\pi}[T_\Delta]<\infty.$$
Recalling that $\eta\ll\lambda$ and using Assumption \ref{Assumption-transition-Q}, it easily follows that $\gamma\ll\lambda$.
Suppose now that $\xi\in\bcal{M}^+(\XX)$ is a solution of \eqref{eq-marginal}. A first direct consequence is that $\xi\ll\lambda$.
Iterating this equation we obtain that
\begin{align}
\xi=\eta \mathbb{I}_{\Delta^c} \sum_{k=0}^{t} (Q_{\pi}\mathbb{I}_{\Delta^c})^k + \xi ( Q_{\pi}\mathbb{I}_{\Delta^c} )^{t+1}\quad\hbox{for any $t\in\NN$}.
\label{expression-marginal-1}
\end{align}
Now,  by Lemma \ref{lemma-Identity-absorbing}(iii)  we have
$\mathbb{I}_{\Delta^c}  (Q_{\pi}\mathbb{I}_{\Delta^c})^k =Q^k_\pi  \mathbb{I}_{\Delta^c}$ and $
( Q_{\pi}\mathbb{I}_{\Delta^c} )^{t+1}=Q^{t+1}_\pi\mathbb{I}_{\Delta^c}$. 
Therefore,  the equation \eqref{expression-marginal-1} becomes
\begin{equation}\label{eq-iterate-t}
\xi=\eta\sum_{k=0}^{t} Q_{\pi}^k\mathbb{I}_{\Delta^c} + \xi  Q_{\pi}^{t+1}\mathbb{I}_{\Delta^c} \quad\hbox{for any $t\in\NN$}.
\end{equation}
Using Lemma \ref{prop-survival-TDelta} (here we make use the fact that $\xi\ll\lambda$) it follows that
$\xi  Q_{\pi}^{t+1}\mathbb{I}_{\Delta^c} (\XX)=\xi  Q_{\pi}^{t+1}(\Delta^c)$ converges to $0$ as $t\rightarrow\infty$. Therefore, taking the limit as  $t\rightarrow\infty$ in equation \eqref{expression-marginal-1} 
we get that indeed $\xi=\gamma$, which completes the proof of the uniqueness. 
For the last statement observe that, by~\eqref{eq-marginal-occupation-measure}, it follows that $\gamma$  is precisely the $\XX$-marginal measure of the occupation measure $\mu_{\eta,\pi}$, that is, $\gamma=\mu^\XX_{\eta,\pi}$.
$\hfill$ $\Box$

\begin{proposition}\label{prop-linear-equation-occupation} 
The occupation measures satisfy the following properties.
\begin{enumerate}
\item[(i).] Given $\pi\in\tilde{\mathbf{\Pi}}$, the occupation measure $\mu_{\eta,\pi}$  satisfies the so-called \textit{characteristic equations} (written in the variable $\mu\in\bcal{M}^+(\XX\times\AA)$):  
\begin{equation}\label{eq-linear-equation}
\mu(\KK^c)=0
\quad
\hbox{and}\quad
\mu^\XX=(\eta+  \mu Q )\mathbb{I}_{\Delta^c}.
\end{equation}
\item[(ii).]  If $\pi\in\tilde{\mathbf{M}}$ is a correlated Markov strategy then $\mu_{\eta,\pi}=\mu^\XX_{\eta,\pi}\otimes \pi$. Moreover, if 
$\pi\in\mathbf{M}$ then $\mu_{\eta,\pi}=\mu_{\eta,\pi}^{\XX\times\AA^{i}}\otimes\pi^{-i}$ for any $i\in\{1,\ldots,N\}$.
\item[(iii).] If  $\pi\in\mathbf{\Pi}$ is such that $\pi^{-i}\in\mathbf{M}^{-i}$ then there exists $\sigma\in\mathbf{M}^i$ with $\mu_{\eta,\pi}=\mu_{\eta,(\pi^{-i},\sigma)}$. 
\item[(iv).] If $\mu\in\bcal{M}^+(\XX\times\AA)$ is a solution of \eqref{eq-linear-equation} then there exists $\pi\in\tilde{\mathbf{M}}$ such that
$\mu=\mu_{\eta,\pi}$ and so,
$$\tilde{\bcal{O}}_{\eta}=  \big\{\mu\in\bcal{M}^+(\XX\times\AA) : \mu(\KK^c)=0 \text{ and } \mu^{\XX}=(\eta+ \mu Q)\mathbb{I}_{\Delta^c} \big\}.$$
Moreover, we have $\mu^{\XX}\ll\lambda$ and if $\lambda\ll\eta$ then $\mu^{\XX}\sim\lambda_{\Delta^c}$.
\end{enumerate}
\end{proposition}
\textbf{Proof.} (i). 
To prove the stated result, note that for any $B\in\bfrak{X}$ we have
\begin{align*}
\mu^\XX(B) = \sum_{t=0}^\infty \mathbb{P}_{\eta,\pi}\{T_\Delta>t,X_t\in B\} 
=  \eta(B-\Delta)+\sum_{t=1}^\infty \mathbb{E}_{\eta,\pi}\big[  \mathbb{P}_{\eta,\pi}\{T_\Delta>t,X_t\in B\mid H_{t-1},A_{t-1}\}\big].
\end{align*}
Observe now that for each $t\ge1$, on the set $\{T_\Delta\le t-1\}$, the conditional probability within brackets vanishes, and so
\begin{eqnarray*}
\mu^\XX(B) 
&= & \eta(B-\Delta)+\sum_{t=1}^\infty \mathbb{E}_{\eta,\pi}\big[  Q(B-\Delta\mid X_{t-1},A_{t-1})
\cdot\mathbf{I}_{\{T_\Delta>t-1\}}\big]\\
&=&  \eta(B-\Delta)+\int_{\XX\times\AA}Q(B-\Delta |x,a)\mu(dx,da),
\end{eqnarray*}
which can be equivalently written precisely as $\mu^\XX=(\eta+\mu Q)\mathbb{I}_{\Delta^c}$.
By construction of the state-action process, it is clear that $\mu(\KK^c)=0$. 
\\[5pt]\noindent
(ii).  Given $B\in\bfrak{X}$ and $D^i\in\bfrak{B}(\AA^i)$ we can write 
\begin{eqnarray*}
\mu_{\eta,\pi}(B\times D^1\times\ldots\times D^N)&=& \sum_{t=0}^\infty  \mathbb{E}_{\eta,\pi} \Big[ \mathbf{I}_{\{T_\Delta>t\}} \mathbf{I}_{\{X_t\in B\}} \pi(D^1\times\ldots\times D^N|X_t)\Big]\\&=&
\int_B  \pi(D^1\times\ldots\times D^N|x)\mu^\XX_{\eta,\pi}(dx)
\end{eqnarray*}
because, precisely, $\mu_{\eta,\pi}^\XX(\Gamma)=\sum_{t\ge0} \mathbb{P}_{\eta,\pi}\{T_\Delta>t,X_t\in \Gamma\}$ for $\Gamma\in\bfrak{X}$,
and the stated result follows. The second part of the statement is an easy consequence of the first part and the fact that, this time, $\pi\in\mathbf{M}$ is a noncooperative Markov profile.
\\[5pt]\noindent
(iii). The occupation measure of the strategy profile $\pi$ satisfies, for $B\in\bfrak{X}$ and $D^j\in\bfrak{B}(\AA^j)$ for $1\le j\le N$
\begin{eqnarray}
\mu_{\eta,\pi}(B\times D^1\times\ldots \times D^N) & =&  \sum_{t=0}^\infty  \mathbb{E}_{\eta,\pi} \Big[ \mathbf{I}_{\{T_\Delta>t\}} \mathbf{I}_{\{X_t\in B\}} \pi^{-i}(D^{-i}|X_t)\pi^i(D^i|H_t)\Big]\nonumber \\
&=& \int_{B} \pi^{-i}(D^{-i}|x) \mu^{\XX\times\AA^{i}}_{\eta,\pi}(dx\times D^{i}),\label{eq-tool-equivalent-markov}
\end{eqnarray}
where  $D^{-i}$ denotes the product of all the sets $D^j$ except $D^i$, and where $\mu^{\XX\times\AA^{i}}_{\eta,\pi}(dx\times D^{i})$ denotes integration with respect to the measure $B\mapsto \mu^{\XX\times\AA^{i}}_{\eta,\pi}(B\times D^{i})$.
By the disintegration result in Lemma \ref{lemma-disintegration}, there exists some $\sigma\in\mathbf{M}^i$ such that    $\mu^{\XX\times\AA^i}_{\eta,\pi}=\mu^{\XX}_{\eta,\pi}\otimes\sigma$. It then follows from \eqref{eq-tool-equivalent-markov} that $\mu_{\eta,\pi}=\mu^\XX_{\eta,\pi}\otimes (\pi^{-i},\sigma)$ and so applying statement (i) in this proposition,
\begin{eqnarray*}
\mu_{\eta,\pi}^\XX&=& (\eta +  \mu_{\eta,\pi} Q)\mathbb{I}_{\Delta^c}\\
&=& (\eta +  \mu^\XX_{\eta,\pi} Q_{(\pi^{-i},\sigma)})\mathbb{I}_{\Delta^c}.
\end{eqnarray*}
By Lemma \ref{prop-linear-equation-occupation-unicity}   we derive that   $\mu_{\eta,\pi}^\XX=\mu^\XX_{\eta,(\pi^{-i},\sigma)}$,  and by item (ii) that $\mu_{\eta,\pi}=\mu_{\eta,(\pi^{-i},\sigma)}$.
\\[5pt]\noindent
(iv). By the disintegration in Lemma \ref{lemma-disintegration}, we obtain that $\mu=\mu^\XX\otimes\pi$ for some $\pi\in\tilde{\mathbf{M}}$. Therefore,
$\mu^\XX$ satisfies equation \eqref{eq-marginal} and so
$\mu^\XX=\mu^\XX_{\eta,\pi}\ll\lambda$, while item (ii) yields 
$\mu=\mu^\XX_{\eta,\pi}\otimes\pi=\mu_{\eta,\pi}$. By using item (i), we get the characterization of $\tilde{\bcal{O}}_{\eta}$.
Now, let $B\in\bfrak{X}$ be  such that $B\subseteq\Delta^c$ and $\mu^\XX(B)=0$. Since \eqref{eq-linear-equation}
implies that $\mu^\XX(B)\ge \eta_{\Delta^c}(B)$ then necessarily
 $\lambda_{\Delta^c}(B)=0$. We conclude that $\mu^\XX\sim \lambda_{\Delta^c}$.
\hfill$\Box$

\bigskip

Now, we introduce $\tilde{\bcal{O}}_\eta^i$ as the set of possible responses for each player $1\le i\le N$.
\begin{definition}
\label{Definition-Otilde-i}
Given an initial distribution $\eta\in\bcal{P}_\lambda(\XX)$ we define  
$$\tilde{\bcal{O}}_\eta^i=\{\mu^{\XX\times\AA^i}:\mu\in\tilde{\bcal{O}}_\eta\}\subseteq\bcal{M}^+(\XX\times\AA^i).$$
\end{definition}

In our next result, we use the notion of a uniformly absorbing game model (see Definition~\ref{Game-absorbing}). Recall that, by Assumption \ref{Assump}, we are considering an initial distribution $\eta\in\bcal{P}_\lambda(\XX)$ such that the game model $\mathcal{G}(\eta,\rho)$ is absorbing to $\Delta$. 
To obtain this important result, we need two preliminary technical Lemmas.
A direct consequence of Assumption \ref{Assump} is the following result.
\begin{lemma}\label{lem-preliminary-0}
 If $v\in L^{\infty}(\XX,\bfrak{X},\lambda)$ then $Qv\in\car_b(\XX\times\AA,\RR)$.
\end{lemma}
In our next lemma, recall that  $\tilde{\bcal{Y}}$ is endowed with the narrow topology   and that  in $ L^\infty(\XX,\bfrak{X},\lambda)$ we consider the weak$^*$ convergence. 
\begin{lemma}\label{lem-trace-2}
The following continuity results hold.
\begin{itemize}
\item[(i).] 
Given any $f\in\car_b(\XX\times\AA)$ and $v\in L^\infty(\XX,\bfrak{X},\lambda)$, the mappings 
$\pi\mapsto f_\pi$ and $\pi\mapsto Q_\pi v$ from  $\tilde{\bcal{Y}}$ to $ L^\infty(\XX,\bfrak{X},\lambda)$
are continuous.
\item[(ii).] If $v_n\wstar v$ in  $ L^\infty(\XX,\bfrak{X},\lambda)$ and $\pi_n\rightarrow\pi$  in  $\tilde{\bcal{Y}}$ then, for any $t\ge0$, we have $Q^t_{\pi_n}v_n\wstar Q_\pi v$.
\end{itemize}
\end{lemma}
\textbf{Proof.} Part (i) is a direct consequence of the definition of the narrow topology and Lemma \ref{lem-preliminary-0}. For item (ii), the reader is referred to Lemma~4.1 in \cite{nowak19} or Lemmas 3.6 and 3.7 in \cite{dufour22}. \hfill$\Box$
\\[10pt]\indent
The above result just concerns convergence in $\tilde{\bcal{Y}}$ and it is not necessarily true for convergence in $\bcal{Y}$. Indeed, the fact that $\bcal{Y}\subseteq\tilde{\bcal{Y}}$ should not be misleading since, in~$\bcal{Y}$, we are considering the product topology of the $\bcal{Y}^i$ which, as  mentioned in Remark \ref{rem-trace}, does not coincide with the trace topology of $\tilde{\bcal{Y}}$. 

The next proposition shows that it is necessary to reinforce the hypothesis of an absorbing model
by assuming that the model is uniformly absorbing: we will need this to show that the set~$\tilde{\bcal{O}}_\eta^i$ of possible responses of each player
is compact in order to use the Kakutani-Fan-Glicksberg fixed point theorem leading to the existence of a Markovian noncooperative equilibrium.

\begin{proposition}
\label{Compactness-set-D} 
The sets $\tilde{\bcal{O}}_\eta$ and $\tilde{\bcal{O}}_\eta^i$  are convex and the following statements are equivalent.
\begin{itemize}
\item[$(a)$] The game model $\mathcal{G}(\eta,\rho)$ is uniformly absorbing to $\Delta$.
\item[$(b)$] The set $\tilde{\bcal{O}}_\eta$ is a compact metric space for the $ws$-topology.
\item[$(c)$] The set $\tilde{\bcal{O}}_\eta^i$ is a compact metric space for the $ws$-topology with $i\in\{1,\ldots,N\}$.
\end{itemize}
\end{proposition}
\textbf{Proof.} 
Regarding the first claim, observe that the convexity of $\tilde{\bcal{O}}_\eta$ is a direct consequence of Proposition \ref{prop-linear-equation-occupation}(iv).
Convexity of $\tilde{\bcal{O}}_\eta^i$ is a straightforward consequence of convexity of $\tilde{\bcal{O}}_\eta$. 

\bigskip

\noindent
$(a)\Rightarrow (b)$
Let us first show that $\tilde{\bcal{O}}_\eta$ is relatively compact for the $ws$-topology.
Applying Theorem 5.2.(ii) in \cite{balder01}, this is equivalent to show that the set of $\XX$-marginal measures of $\tilde{\bcal{O}}_\eta$, which we denote by $\tilde{\bcal{O}}_\eta^\XX$, is relatively $s$-compact and that the set of $\AA$-marginal measures of $\tilde{\bcal{O}}_\eta$, denoted by 
$\tilde{\bcal{O}}_\eta^\AA$, is relatively $w$-compact. 
Recalling Remark \ref{rem-occupation-measures}(a), we have that $\tilde{\bcal{O}}_\eta$ is a bounded subset of $\bcal{M}^+(\XX\times\AA)$.
Since~$\AA$ is compact, it is clear that $\tilde{\bcal{O}}_\eta^\AA$ is relatively $w$-compact by using   \cite[Theorem~8.6.7]{bogachev07}.
To prove that~$\tilde{\bcal{O}}_\eta^\XX$ is relatively $s$-compact, let us show that 
\begin{equation}\label{eq-to-show-relatively-compact}
\lim_{n\rightarrow\infty} \sup_{\mu\in\tilde{\bcal{O}}_\eta} \mu^\XX(\Gamma_{n})=0
\end{equation}
for any decreasing sequence of sets $\Gamma_{n}\in\bfrak{X}$ with $\Gamma_{n}\downarrow\emptyset$.
Indeed, from \cite[Lemma 4.6.5]{bogachev07} this implies that $\tilde{\bcal{O}}_\eta^\XX$ is uniformly countably additive and so,
relatively compact for the $s$-topology; see \cite[Theorem 4.7.25]{bogachev07}.
Since $\mu^\XX(\Delta)=0$, there is no loss of generality in assuming that the  $\Gamma_n$ are subsets of $\Delta^c$.  By Proposition \ref{prop-linear-equation-occupation}(iv),  for every $\mu\in\tilde{\bcal{O}}_\eta$ there exists a correlated Markov strategy $\pi_\mu\in\tilde{\bcal{Y}}$  such that $\mu=\mu_{\eta,\pi_\mu}$ and so, 
for any fixed $k\ge0$,
\begin{eqnarray*}
\mu^\XX(\Gamma_n)  &\le&  \sum_{t=0}^k  \mathbb{P}_{\eta,\pi_\mu}\{X_t\in\Gamma_n\}+\sum_{t>k}\mathbb{P}_{\eta,\pi_\mu}\{X_t\in\Delta^c\}\\
&=&  \sum_{t=0}^k \mathbb{P}_{\eta,\pi_\mu}\{X_t\in\Gamma_n\}+\sum_{t>k}\mathbb{P}_{\eta,\pi_\mu}\{T_\Delta>t\}
\end{eqnarray*}
Therefore,
\begin{equation}\label{eq-tool-relatively-compact}
\sup_{\mu\in\tilde{\bcal{O}}_\eta} \mu^\XX(\Gamma_n)\le 
\sum_{t=0}^k \sup_{\mu\in\tilde{\bcal{O}}_\eta}  \mathbb{P}_{\eta,\pi_\mu}\{X_t\in\Gamma_n\}  +
\sup_{\mu\in\tilde{\bcal{O}}_\eta} \sum_{t>k}\mathbb{P}_{\eta,\pi_\mu}\{T_\Delta>t\}.
\end{equation}
Let us now pay attention to first term in righthand of \eqref{eq-tool-relatively-compact}. Suppose first that $0\le t\le k$ and $n\in\NN$ remain fixed. By Lemma \ref{lem-trace-2} we have that the mapping $\pi\mapsto Q^t_\pi(\Gamma_n|\cdot)$ from $\tilde{\bcal{Y}}$ to $L^\infty(\XX,\bfrak{X},\lambda)$ is continuous. Since $\eta\ll\lambda$, this implies that the mapping $\pi\mapsto \mathbb{P}_{\eta,\pi_\mu}\{X_t\in\Gamma_n\}$ is continuous on~$\tilde{\bcal{Y}}$.
Observe now that, by hypothesis, we have
 $\mathbf{I}_{\Gamma_n}\wstar0$ as $n\rightarrow\infty$, and so by Lemma~\ref{lem-trace-2} again, for every $\pi\in\tilde{\bcal{Y}}$ we have $Q^t_\pi(\Gamma_n|\cdot)\wstar0$ and, therefore, $\mathbb{P}_{\eta,\pi}\{X_t\in\Gamma_n\}\rightarrow0$.
 Summarizing, the sequence (in $n\in\NN$) of continuous mappings $\pi\mapsto \mathbb{P}_{\eta,\pi}\{X_t\in\Gamma_n\}$ decreases to $0$ and hence, by Dini's theorem, the convergence is uniform since $\tilde{\bcal{Y}}$ is a compact metric space. So, for each fixed $0\le t\le k$ we have
$$\lim_{n\rightarrow\infty}  \sup_{\mu\in\tilde{\bcal{O}}_\eta}  \mathbb{P}_{\eta,\pi_\mu}\{X_t\in\Gamma_n\} =0.$$
Regarding the rightmost expression in \eqref{eq-tool-relatively-compact}, it converges to $0$ as $n\rightarrow\infty$ as a direct consequence of the fact
that $\mathcal{G}(\eta,\rho)$ is uniformly absorbing to $\Delta$. This completes the proof of \eqref{eq-to-show-relatively-compact}.
Therefore, once we know that $\tilde{\bcal{O}}_\eta$ is relatively compact for the $ws$-topology, it follows that it is also metrizable by Proposition 2.3 in \cite{balder01}.

To prove the compactness of $\tilde{\bcal{O}}_\eta$, the last step consists in showing that it is closed. To see this, let $\{\mu_n\}_{n\ge0}$ be a sequence in $\tilde{\bcal{O}}_\eta$ converging in the $ws$-topology to some 
 $\mu\in\bcal{M}^+(\XX\times\AA)$. 
 First of all, let us show that $\mu(\KK^c)=0$. The measurable function $(x,a)\mapsto \mathbf{I}_{\KK^c}(x,a)$  is such that
$a\mapsto \mathbf{I}_{\KK^c}(x,a)=\mathbf{I}_{\AA^c(x)}(a)$
is lower semicontinuous on $\AA$ for any fixed $x\in\XX$ because $\AA(x)$ is compact. 
Thus,
$\mathbf{I}_{\KK^c}$ is a nonnegative normal integrand  and \cite[Theorem 3.1.(c)]{balder01} yields
$\liminf_{n} \mu_{n}(\KK^c)\geq \mu(\KK^c)$ and so  $\mu(\KK^c)=0$. On the other hand, it is clear that $\mu_n^\XX(\Delta)=0$ for all $n\ge0$ implies that $\mu^\XX(\Delta)=0$. To conclude the proof, choose an arbitrary measurable subset $B$ of $\Delta^c$. For every $n\ge0$ we have
$$\mu_n^\XX(B)=\eta(B)+ \int_{\XX\times\AA} Q\mathbf{I}_B(x,a)  \mu_n(dx,da)
$$
By Lemma~\ref{lem-preliminary-0}, the function $Q\mathbf{I}_B(x,a)$  is in  $\car_b(\XX\times\AA,\RR)$ so that we can take limits as $n\rightarrow\infty$ to obtain that $\mu^\XX(B)=\eta(B)+\mu Q(B)$, thus completing the proof that $\mu\in\tilde{\bcal{O}}_\eta$.

\bigskip

\noindent
$(b)\Rightarrow (c)$
Since the mapping from $\bcal{M}^+(\XX\times\AA)$ to $\bcal{M}^+(\XX\times\AA^{i})$ which associates to  $\mu\in\bcal{M}^+(\XX\times\AA)$ its marginal measure $\mu^{\XX\times\AA^i}$ is continuous for the respective $ws$-topologies, it follows that   $\tilde{\bcal{O}}_\eta^i$ is compact.  Again from \cite[Proposition 2.3]{balder01}, noting that the set of $\XX$-marginal  measures of $\tilde{\bcal{O}}_\eta^i$ is precisely $\tilde{\bcal{O}}^\XX_{\eta}$, which has  been shown to be relatively $s$-compact, we conclude that $\tilde{\bcal{O}}_\eta^i$ is metrizable. 

\bigskip

\noindent
$(b)\Rightarrow (a)$
Since $\tilde{\bcal{O}}_\eta$ is compact for the $ws$-topology, it follows from Theorem 5.2 in \cite{balder01} that the set of $\XX$-marginal measures of $\tilde{\bcal{O}}_\eta$ (denoted by $\tilde{\bcal{O}}_\eta^\XX$)
is relatively s-compact.
By Proposition \ref{prop-linear-equation-occupation}(iv), $\tilde{\bcal{O}}_\eta^\XX=\{\mu^\XX_{\pi} : \pi\in \tilde{\bcal{Y}}\}$.
Combining Proposition 2.2 in \cite{balder01} and Corollary 2.7 in \cite{ganssler71} we get that the family $\{h_\pi\}_{\pi\in\tilde{\bcal{Y}}}$ of density functions 
$h_\pi={d\mu^\XX_\pi}/{d\lambda}$ is uniformly $\lambda$-integrable.
Now, observe that for $\pi\in\tilde{\bcal{Y}}$,
$$ \sum_{k=t}^{\infty} \mathbb{P}_{\eta,\pi}\{T_\Delta>k\}= \mu^\XX_{\pi} Q_{\pi}^{t}(\Delta^c)=\int_{\XX} Q_{\pi}^{t}(\Delta^c |x) h_\pi(x) \lambda(dx)$$
and by using Lemma \ref{prop-survival-TDelta} we can conclude that the rightmost term in the previous equation converges to zero uniformly in $\pi\in\tilde{\bcal{Y}}$ as $t\rightarrow\infty$.
This establishes that $\mathcal{G}(\eta,\rho)$ is indeed uniformly absorbing to $\Delta$.

\bigskip

\noindent
$(c)\Rightarrow (a)$
Observe that for $i\in\{1,\ldots,N\}$, the sets of $\XX$-marginal measures of $\tilde{\bcal{O}}_\eta$ and $\tilde{\bcal{O}}_\eta^i$ are the same.
Consequently, by Theorem 5.2 in \cite{balder01} the set of $\XX$-marginal measures of $\tilde{\bcal{O}}_\eta$ is relatively s-compact and the rest of the proof is identical to that of $(b)\Rightarrow (a)$.
\hfill$\Box$

\subsection{Continuity properties of the occupation measures}
\textit{In this subsection we shall assume that the game model $\mathcal{G}(\eta,\rho)$ satisfies Assumptions \ref{Assump-2} and~\ref{Assumption-B}.}
In particular under condition \ref{Assump-absorbing-uniformly}, the set $\tilde{\bcal{O}}_\eta$ of occupation measures is compact by Proposition~\ref{Compactness-set-D}.
Under condition \ref{Assumption-B}, we can obtain the following result similar to Lemma \ref{lem-trace-2}(ii), whose proof is omitted.
\begin{lemma}\label{lem-preliminary-bis}
If $v_n\wstar v$  in $L^{\infty}(\XX,\bfrak{X},\lambda)$ and $\pi_n\rightarrow\pi$ in $\bcal{Y}$ then
$Q^t_{\pi_n} v_{n} \wstar Q^t_{\pi}v$ for any $t\in\NN$.
\end{lemma}
At this point, recall  the notation $\eta_n$ (see~\eqref{eq-def-zeta-n}) for the initial distributions which are a combination of $\eta$ and $\lambda$. 
\begin{proposition} 
\label{Convergence-mu-general} Under any of the conditions (i) and (ii) below: 
\begin{itemize}
\item[(i).] $\pi_n\rightarrow\pi$ in $\bcal{Y}$ and $f\in\car_b(\XX\times\AA,\RR)$ is such that $\mathbb{I}_{\Delta^{c}}f_{\pi_{n}}(\cdot)\wstar \mathbb{I}_{\Delta^{c}}f_{\pi}(\cdot)$ in $L^{\infty}(\XX,\bfrak{X},\lambda)$,
\item[(ii).] $\pi_n\rightarrow\pi$ in $\tilde{\bcal{Y}}$ and $f\in\car_b(\XX\times\AA,\RR)$,
\end{itemize}
we have the following limits:
$$
\lim_{n\rightarrow\infty} \int_{\XX\times\AA} fd\mu_{\eta_n,\pi_{n}} = \int_{\XX\times\AA} fd\mu_{\eta,\pi}
\quad\hbox{and}\quad  
\lim_{n\rightarrow\infty} \int_{\XX\times\AA} fd\mu_{\eta,\pi_{n}} = \int_{\XX\times\AA} fd\mu_{\eta,\pi}
$$
\end{proposition}
\textbf{Proof:}
We will only prove the first limit in case (i), the remaining cases being obtained by using similar arguments.
Recalling \eqref{eq-def-zeta-n}, observe that $$\mu_{\eta_n,\pi_{n}}=\frac{n}{n+1} \mu_{\eta,\pi_{n}} + \frac{1}{n+1} \mu_{\lambda,\pi_{n}}\quad\hbox{and}\quad \int_{\XX\times\AA} f(x,a)d\mu_{\lambda,\pi_{n}}\leq \mathbf{f} \sup_{\pi\in\tilde{\mathbf{\Pi}}} \mu_{\lambda,\pi}(\XX\times\AA)$$ for some constant $\mathbf{f}$.
From Remark \ref{rem-occupation-measures}(a), we only have to show that 
\begin{equation}\label{eq-to-show}
 \lim_{n\rightarrow\infty} \int_{\XX\times\AA} f(x,a)d\mu_{\eta,\pi_{n}}= \int_{\XX\times\AA} f(x,a)d\mu_{\eta,\pi}.
 \end{equation}
Equivalently, the above sequence being bounded, we will prove that any convergent subsequence 
has the desired limit. 
To simplify the notation, and without loss of generality, we will suppose that the whole sequence is converging and also that 
$\{\mu_{\eta,\pi_{n}}\}_{n\in\NN}$ is a convergent sequence in $\tilde{\bcal{O}}_\eta$  (recall Assumption \ref{Assump-absorbing-uniformly} and Proposition \ref{Compactness-set-D}).
We have $\mu_{\eta,\pi_{n}}=\mu^{\mathbf{X}}_{\eta,\pi_{n}}\otimes \pi_{n}$ with (by Lemma \ref{prop-linear-equation-occupation-unicity} and, in particular, \eqref{eq-iterate-t})
$$\mu^{\mathbf{X}}_{\eta,\pi_{n}}=\sum_{k=0}^{t-1} \eta Q^k_{\pi_{n}}\mathbb{I}_{\Delta^{c}} + \mu^{\mathbf{X}}_{\eta,\pi_{n}}Q^t_{\pi_{n}}\mathbb{I}_{\Delta^{c}}$$
for any $t\in\NN^{*}$. Consequently, integrating the function $f_{\pi_n}$ with respect to the above measures, we can write 
\begin{eqnarray}
\int_{\XX\times\AA} f d\mu_{\eta,\pi_{n}}&=& 
\sum_{k=0}^{t-1} \int_{\XX}f_{\pi_{n}}\, d  \eta Q^k_{\pi_{n}}\mathbb{I}_{\Delta^{c}} 
+\int_{\XX}  f_{\pi_n}\, d \mu^{\mathbf{X}}_{\eta,\pi_{n}} Q^t_{\pi_{n}}\mathbb{I}_{\Delta^{c}}  \nonumber
\\
&=& 
\sum_{k=0}^{t-1} \int_{\XX} Q^k_{\pi_{n}}\mathbb{I}_{\Delta^{c}}f_{\pi_{n}}d  \eta 
+\int_{\XX} \mathbb{I}_{\Delta^{c}}f_{\pi_{n}} d \mu^{\mathbf{X}}_{\eta,\pi_{n}}Q^t_{\pi_{n}},
\label{Convergence-mu-eq1}
\end{eqnarray}
for any $t\in\NN^*$.
Observe also that 
\begin{eqnarray}
\Big|  \int_{\XX} \mathbb{I}_{\Delta^{c}}f_{\pi_{n}}\, d \mu^{\mathbf{X}}_{\eta,\pi_{n}}Q^t_{\pi_{n}} \Big| \leq \mathbf{f}\cdot \mu^{\mathbf{X}}_{\eta,\pi_{n}}Q^t_{\pi_{n}}(\Delta^{c})
= \mathbf{f} \sum_{k=t}^{\infty} \mathbb{P}_{\eta,\pi_{n}} \{T_\Delta>k\}
\label{Convergence-mu-eq2}
\end{eqnarray}
for any $t\in\NN^*$ and $n \in\NN$.
By hypothesis we have  $\mathbb{I}_{\Delta^{c}}f_{\pi_{n}}(\cdot)\wstar \mathbb{I}_{\Delta^{c}}f_{\pi}(\cdot)$
in $L^{\infty}(\XX,\bfrak{X},\lambda)$
 and since $d\eta/d\lambda$ is in $L^{1}(\XX,\bfrak{X},\lambda)$, we have by Lemma \ref{lem-preliminary-bis} that
\begin{eqnarray}
\int_{\XX} Q^k_{\pi_{n}} \mathbb{I}_{\Delta^{c}}f_{\pi_{n}}\, \frac{d\eta}{d\lambda}\, d \lambda \rightarrow \int_{\XX} Q^k_{\pi} \mathbb{I}_{\Delta^{c}}f_{\pi}\, \frac{d\eta}{d\lambda}\, d \lambda.
\label{Convergence-mu-eq3}
\end{eqnarray}
Combining equations \eqref{Convergence-mu-eq1}--\eqref{Convergence-mu-eq3} we get that for any $t\in\NN^*$,
\begin{eqnarray*}
\Big| \lim_{n\rightarrow\infty} \int_{\XX\times\AA} fd\mu_{\eta,\pi_{n}} -\sum_{k=0}^{t-1} \int_{\XX} f_\pi \, d\eta Q^k_{\pi}\mathbb{I}_{\Delta^{c}} \Big|
\leq \mathbf{f} \sup_{n\in\NN}   \sum_{k=t}^{\infty} \mathbb{P}_{\eta,\pi_{n}} \{T_\Delta>k\}.
\end{eqnarray*}
Finally, by Assumption \ref{Assumption-B} we get
$$\lim_{n\rightarrow\infty} \int_{\XX\times\AA} fd\mu_{\eta,\pi_{n}} = \sum_{k=0}^{\infty} \int_{\XX} f_\pi \, d\eta Q^k_{\pi}\mathbb{I}_{\Delta^{c}}
=\int_{\XX} f_\pi d\mu^\XX_{\eta,\pi}
=\int_{\XX\times\AA} fd\mu_{\eta,\pi},$$  this establishes the limit in \eqref{eq-to-show}.
\hfill $\Box$
\\[10pt]\indent
The next result will be useful in the forthcoming.
\begin{corollary}
\label{Convergence-mu}
The following convergence results hold.
\begin{enumerate}
 \item[(i).] If  $\{\pi_n\}$ in $\tilde{\bcal{Y}}$ converges to $\pi\in\tilde{\bcal{Y}}$ then $\mu_{\eta,\pi_n}\rightarrow \mu_{\eta,\pi}$ and $\mu_{\eta_n,\pi_n}\rightarrow \mu_{\eta,\pi}$.
\item[(ii).] Consider $\pi^{-i}\in\bcal{Y}^{-i}$ fixed for $i\in\{1,\ldots,N\}$. If $\sigma_n\rightarrow\sigma$ in $\bcal{Y}^i$ then $\mu_{\eta,(\pi^{-i},\sigma_n)}\rightarrow \mu_{\eta,(\pi^{-i},\sigma)}$ and  $\mu_{\eta_n,(\pi^{-i},\sigma_n)}\rightarrow \mu_{\eta,(\pi^{-i},\sigma)}$
 in the $ws$-topology.
 \item[(iii).] Fix a player $1\le i\le N$ and suppose that $\{\pi_n\}$ in $\bcal{Y}$ converges to $\pi\in\bcal{Y}$. 
Then,
$$\mu^{\XX\times\AA^i}_{\eta,\pi_n}\rightarrow \mu_{\eta,\pi}^{\XX\times\AA^i}\quad\hbox{and}\quad
\mu^{\XX\times\AA^i}_{\eta_n,\pi_n}\rightarrow \mu_{\eta,\pi}^{\XX\times\AA^i}.
$$
 \item[(iv).] Consider the initial distributions $\eta_n$ defined in \eqref{eq-def-zeta-n} and an arbitrary sequence $\pi_n\rightarrow\pi$ in $\bcal{Y}$.
Then,
$$R^i(\eta,\pi_{n})\rightarrow R^i(\eta,\pi)\quad\hbox{and}\quad
C^i(\eta,\pi_{n})\rightarrow C^i(\eta,\pi).
$$
and also
$$R^i(\eta_n,\pi_{n})\rightarrow R^i(\eta,\pi)\quad\hbox{and}\quad
C^i(\eta_n,\pi_{n})\rightarrow C^i(\eta,\pi).
$$
 \end{enumerate}
\end{corollary}
\textbf{Proof:}
Part (i) follows directly from the second condition in Proposition \ref{Convergence-mu-general}.
For the second part, observe that  if $g\in\car(\XX\times\AA)$ then
$$g^i(x,a^i)=\int_{\AA^{-i}} g(x,(a^i,a^{-i}))\pi^{-i}(da^{-i}|x)$$ is in $\car(\XX\times\AA^i)$.
With this in mind, it is easy to see that, letting
 $\pi_n=(\pi^{-i},\sigma_n)$ and $\pi=(\pi^{-i},\sigma)$, we have that $\pi_n\rightarrow\pi$ in~$\tilde{\bcal{Y}}$. Apply now Proposition \ref{Convergence-mu-general}(ii).
For item (iii), consider a fixed integer $i\in\{1,\ldots, N\}$ and an arbitrary $g\in \car_b(\XX\times\AA^i,\RR)$.
Define then the  function $f$ on  $\XX\times\AA$ by $f(x,a_{1},\ldots,a_{i},\ldots,a_{N})=g(x,a_{i})$ with $f\in\car_b(\XX\times\AA,\RR)$, which satisfies
 $\mathbb{I}_{\Delta^{c}}f_{\pi_{n}}(\cdot)\wstar \mathbb{I}_{\Delta^{c}}f_{\pi}(\cdot)$ in $L^{\infty}(\XX,\bfrak{X},\lambda)$. Applying
Proposition \ref{Convergence-mu-general} we conclude that 
$\int_{\XX\times\AA^i} g d\mu^{\XX\times\AA^i}_{\eta,\pi_n} \rightarrow \int_{\XX\times\AA^i} g d\mu^{\XX\times\AA^i}_{\eta,\pi}$,
and the result follows. 
The proof of item (iv) is a direct consequence of Assumption \ref{Assumption-B}, Remark \ref{rem-occupation-measures}(c),  and Proposition~\ref{Convergence-mu-general}.
\hfill $\Box$
\\[10pt]\indent
By Corollary \ref{Convergence-mu}(i), when considering the narrow convergence in $\tilde{\bcal{Y}}$ we have that the mapping which associates to each $\pi\in\tilde{\bcal{Y}}$ its occupation measure  $\mu_{\eta,\pi}$ is continuous. For the product narrow topology on $\bcal{Y}$, however, such a result is not true in general and, instead, we get weaker results as in Corollary \ref{Convergence-mu}(ii) and (iii) above. Notice that  there is some kind of duality in these  results: indeed, item (ii) shows that if the convergence $\pi_n\rightarrow\pi$ takes place in \textit{just one} variable then we get convergence of the \textit{whole} occupation measures $\mu_{\eta,\pi_n}\rightarrow\mu_{\eta,\pi}$, while in item (iii)  if the \textit{whole} sequence $\pi_n$ converges to $\pi$  then we get convergence of, \textit{individually}, each component of the occupation measures $\mu^{\XX\times\AA^i}_{\eta,\pi_n}\rightarrow \mu^{\XX\times\AA^i}_{\eta,\pi}$.

\section{Proofs of the main results}\label{section4}
\setcounter{equation}{0}

We will show Proposition \ref{th-main-0} and Theorem \ref{th-main} in the constrained case. The unconstrained case is easily deduced from the constrained one.
Indeed, by considering constraint constants satisfying
$\rho < -\mathbf{r} \sup_{\pi\in\Pi}\EE_{\eta,\pi}[T_\Delta] \1$ we get that $C^i(\eta,\pi) > \rho$ for any $\pi\in\Pi$ and any player $i$
(see item (a) of Remark \ref{rem-after-assumption-A}) yielding that the constraints and and the Slater condition \ref{Assump-Slater} are trivially satisfied.

\subsection{Proof of Proposition \ref{th-main-0}}
\label{Proof-Main1}
\textit{We will suppose in this subsection that we are given an initial distribution $\eta\in\bcal{P}(\XX)$ satisfying $\lambda\ll \eta$ and constraint constants $\rho\in\RR^{Np}$ such that the game $\mathcal{G}(\eta,\rho)$ satisfies Assumptions~\ref{Assump-2} and \ref{Assumption-B}. }
By recalling Assumption \ref{Assump-Abs-Cont} and Proposition \ref{prop-linear-equation-occupation}(iv), this yields to an important property, namely, $\mu^{\XX}_{\eta,\pi}\sim\lambda_{\Delta^c}$
for any $\pi\in\tilde{\mathbf{M}}$.
The objective of this subsection is to introduce a correspondence defined as the composition of a function
$\mathcal{J}_\eta:\tilde{\bcal{O}}_\eta^1\times\ldots\times\tilde{\bcal{O}}^N_\eta\rightarrow \bcal{Y}^1\times\ldots\times\bcal{Y}^N$
and a correspondence $\mathcal{H}_{\eta,\rho}: \bcal{Y}^1\times\ldots\times\bcal{Y}^N\tto \tilde{\bcal{O}}_\eta^1\times\ldots\times\tilde{\bcal{O}}^N_\eta$
and to show that it has a fixed point, from which we will derive equilibrium stationary Markov policies in the special case where $\lambda\ll \eta$.

\paragraph{The function $\mathcal{J}_{\eta}$.}
Consider a fixed $\vartheta\in{\mathbf{M}}$.
For any $\pi\in\tilde{\mathbf{M}}$, let us define $\gamma_{\pi}\in\tilde{\mathbf{M}}$ as
\begin{eqnarray}
\gamma_{\pi}(B|x)=\pi(B|x) \mathbf{I}_{\Delta^{c}}(x)+\vartheta(B|x)\mathbf{I}_{\Delta}(x)\quad\hbox{for $B\in\bfrak{B}(\AA)$ and $x\in\XX$.}
\label{Def-Gamma-pi}
\end{eqnarray}
This definition ensures that $\gamma_\pi\in\mathbf{M}$ if $\pi\in\mathbf{M}$.

\begin{lemma}
\label{lemma-preliminary-1} 
\begin{enumerate}
\item[(i).] Given any $\pi\in\tilde{\mathbf{M}}$ we have  equality of the occupation measures $\mu_{\eta,\pi}=\mu_{\eta,\gamma_{\pi}}$.
\item[(ii).]
Let $i\in\{1,\ldots,N\}$ be fixed. 
For any  $m\in\tilde{\bcal{O}}_\eta^i$ consider  
$\tilde{\mathbf{M}}_{m}^i=\{\pi\in\tilde{\mathbf{M}}: \mu^{\XX\times\AA^{i}}_{\eta,\pi}=m\}$.
Then the set $$\{\gamma_\pi^{\AA^i}:\pi\in\tilde{\mathbf{M}}_{m}^i\}\subseteq\mathbf{M}^i$$ is contained in a unique class of equivalence of $\bcal{Y}^i$, that will be denoted by $\mathcal{J}^i_{\eta}(m)$.
\item[(iii).] Given $\pi\in{\mathbf{M}}$ and $1\le i\le N$, let $\sigma=\gamma_{\pi}^{\AA^i}$ and consider $\pi'=(\pi^{-i},\sigma)$. Then $\mu_{\eta,\pi}=\mu_{\eta,\pi'}$.
\end{enumerate}
\end{lemma}
\textbf{Proof.} 
(i). By its definition, it is clear that $Q_\pi\mathbb{I}_{\Delta^c}=Q_{\gamma_\pi}\mathbb{I}_{\Delta^c}$. Combining Lemma \ref{lemma-Identity-absorbing}(iii) and
Lemma \ref{prop-linear-equation-occupation-unicity} we conclude that the $\XX$-marginals of the occupation measures of $\pi$ and $\gamma_\pi$ coincide: $\mu_{\eta,\pi}^\XX=\mu_{\eta,\gamma_{\pi}}^\XX$.  Since $\mu^\XX_{\eta,\gamma_\pi}(\Delta)=0$ and $\pi(\cdot|x)=\gamma_\pi(\cdot|x)$ when $x\in\Delta^c$, we conclude that 
$\mu_{\eta,\gamma_\pi}\otimes\pi=\mu_{\eta,\gamma_\pi}\otimes\gamma_\pi$. Summarizing, we have
$$\mu_{\eta,\pi}=\mu_{\eta,\pi}^\XX\otimes\pi= \mu_{\eta,\gamma_\pi}^\XX\otimes\pi=
 \mu_{\eta,\gamma_\pi}^\XX\otimes\gamma_\pi=\mu_{\eta,\gamma_\pi}$$
 by using Proposition \ref{prop-linear-equation-occupation}(ii).

(ii).
We must show that if $\pi,\bar\pi\in \tilde{\mathbf{M}}$ are such that 
$\mu^{\XX\times\AA^{i}}_{\eta,{\pi}}=\mu^{\XX\times\AA^{i}}_{\eta,\bar\pi}$ then $\gamma_{\pi}^{\AA^{i}}$ and $\gamma_{{\bar\pi}}^{\AA^{i}}$ belong to the same class of equivalence in $\bcal{Y}^i$.
Clearly, we have $\mu^{\XX}_{\eta,{\pi}}=\mu^{\XX}_{\eta,\bar\pi}$.
By Proposition \ref{prop-linear-equation-occupation}(ii), it follows that
$\mu^{\XX\times\AA^{i}}_{\eta,{\pi}}=\mu^{\XX}_{\eta,{\pi}} \otimes \pi^{\AA^{i}}$ and $\mu^{\XX\times\AA^{i}}_{\eta,{\bar \pi}}=\mu^{\XX}_{\eta,{\bar \pi}} \otimes \bar{\pi}^{\AA^{i}}$
implying that $\mu^{\XX}_{\eta,{\pi}} \otimes \pi^{\AA^{i}}=\mu^{\XX}_{\eta,{\pi}} \otimes \bar{\pi}^{\AA^{i}}$.
Now, recalling that $\mu^{\XX}_{\eta,\pi}\sim\lambda_{\Delta^c}$ we obtain by the disintegration Lemma that 
$\pi^{\AA^{i}}(\cdot |x)=\bar{\pi}^{\AA^{i}}(\cdot | x)$ for $\lambda_{\Delta^c}$-almost all $x\in\XX$
which shows that $\gamma^{\AA^{i}}_\pi(\cdot|x)=\gamma^{\AA^{i}}_{\bar{\pi}}(\cdot|x)$  with $\lambda$-almost all $x\in\XX$.

(iii). Proceeding as in part (i), we can show that $\mu_{\eta,\pi}^\XX=\mu_{\eta,\pi'}^\XX$, which is a measure equivalent to $\lambda_{\Delta^c}$. Since $\mu_{\eta,\pi}^{\XX}(\Delta)=\mu_{\eta,\pi'}^{\XX}(\Delta)=0$ and $\pi^{\AA^i}=\sigma$ on~$\Delta^c$ we conclude by Proposition \ref{prop-linear-equation-occupation}(ii) that $\mu_{\eta,\pi}^{\XX\times\AA^i}=\mu_{\eta,\pi'}^{\XX\times\AA^i}$and so,
$$\mu_{\eta,\pi}=\mu_{\eta,\pi}^{\XX\times\AA^i}\otimes\pi^{-i}\quad\hbox{and}\quad \mu_{\eta,\pi'}=\mu_{\eta,\pi'}^{\XX\times\AA^i}\otimes\pi^{-i},$$
and the result follows.
\hfill$\Box$
\\[10pt]\indent
Using the result in Lemma \ref{lemma-preliminary-1}(ii), we indeed have defined a function $\mathcal{J}^i_\eta$ from $\tilde{\bcal{O}}^i_\eta$ to $\bcal{Y}^i$. We can therefore consider the function
$$\mathcal{J}_\eta:\tilde{\bcal{O}}_\eta^1\times\ldots\times\tilde{\bcal{O}}^N_\eta\rightarrow \bcal{Y}^1\times\ldots\times\bcal{Y}^N=\bcal{Y}$$
whose components are the $\mathcal{J}_\eta^i$. Based on this lemma and using Remark \ref{rem-occupation-measures}(d), without risk of confusion we can assume that $\gamma_\pi\in\tilde{\bcal{Y}}$ and $\gamma_\pi^{\AA^i}\in\bcal{Y}^i$ for $\pi\in\tilde{\bcal{Y}}$.  Indeed, two Markov correlated strategies $\pi,\pi'$ in the same class of equivalence of $\bcal{Y}$ have the same occupation measure and they yield the same class of equivalence for $\gamma_\pi^{\AA^i}$ and $\gamma_{\pi'}^{\AA^i}$

\begin{proposition}
\label{Correspondence-J}
The function $\mathcal{J}_{\eta}$ is continuous.
\end{proposition}
\textbf{Proof:} We make the proof of the continuity for the function $\mathcal{J}^i_\eta$ for any fixed $1\le i\le N$. Suppose that 
$\{m_n\}_{n\ge0}$ is a sequence in $\tilde{\bcal{O}}_\eta^i$ converging in the $ws$-topology to some $m\in\tilde{\bcal{O}}_\eta^i$.
There exist  $\pi_n$ and $\pi $ in $\tilde{\bcal{Y}}$ such that $m_{n}=\mu^{\XX\times\AA_{i}}_{\eta,\pi_{n}}$ for any $n\in\NN$ and $m=\mu^{\XX\times\AA_{i}}_{\eta,\pi}$.
Our goal is to prove that $\gamma_{\pi_n}^{\AA_{i}}\rightarrow\gamma_{\pi}^{\AA_{i}}$ in $\bcal{Y}^{i}$. Since $\bcal{Y}^i$ is compact, it suffices to show that
this limit holds for any convergent subsequence of  $\{\gamma_{\pi_n}^{\AA_{i}}\}$ (still denoted by $\{\gamma_{\pi_n}^{\AA_{i}}\}$).
There is no loss of generality in assuming that $\pi_n\rightarrow\pi^{*}$ for some
 $\pi^{*}\in \tilde{\bcal{Y}}$ and, therefore, as a direct consequence of the definition in  \eqref{Def-Gamma-pi} we also have  $\gamma_{\pi_n}\rightarrow\gamma_{\pi^{*}}$ and $\gamma^{\AA^i}_{\pi_n}\rightarrow\gamma^{\AA^i}_{\pi^{*}}$.
By Corollary \ref{Convergence-mu}(i) we obtain that
 $$
\mu_{\eta,\pi_{n}}\rightarrow \mu_{\eta,\pi^{*}}
  \quad\hbox{and so}\quad
m_n=\mu^{\XX\times\AA^{i}}_{\eta,\pi_{n}}\rightarrow \mu^{\XX\times\AA^{i}}_{\eta,\pi^{*}}.$$
We deduce that 
 $\mu^{\XX\times\AA_{i}}_{\eta,\pi^{*}}=\mu^{\XX\times\AA_{i}}_{\eta,\pi}=m$ and from Lemma \ref{lemma-preliminary-1}(ii) we conclude that $\gamma^{\AA^i}_{\pi^*}=\gamma^{\AA^i}_{\pi}$, which completes the proof.
\hfill $\Box$

\paragraph{The correspondence $\mathcal{H}_{\eta,\rho}$.}

Fix a player $i$ and an arbitrary $\pi^{-i}\in\bcal{Y}^{-i}$.
Define
$$
\mathcal{L}_{\eta}^i(\pi^{-i}) =  \big\{
 \mu_{\eta,(\pi^{-i},\sigma)}^{\XX\times\AA^{i}} :\sigma\in\bcal{Y}^{i} \big\}
\subseteq\tilde{\bcal{O}}_\eta^i\subseteq\bcal{M}^+(\XX\times\AA^i),
$$
which is the set of $(\XX\times\AA^i)$-marginals of the occupation measures for the initial distribution $\eta$ and the strategy profiles $(\pi^{-i},\sigma)$ as  the policy $\pi^{-i}$ of all the players (but $i$) remain fixed and the policy of player $i$ varies in $\bcal{Y}^i$.  

\begin{proposition}\label{prop-L1-convex-compact}
Given any $1\le i\le N$ and $\pi^{-i}\in\bcal{Y}^{-i}$, 
 the set $\mathcal{L}_{\eta}(\pi^{-i})$ is convex and compact for the $ws$-topology.
\end{proposition}
\textbf{Proof.} Let $\gamma,\gamma'$ in $\mathcal{L}_{\eta}^i(\pi^{-i})$ and fix some $0\le\alpha\le1$. We want to prove that $\alpha\gamma+(1-\alpha)\gamma'\in\mathcal{L}_{\eta}^i(\pi^{-i})$. 
There exist $\sigma,\sigma'\in\bcal{Y}_{i}$
satisfying $\gamma= \mu_{\eta,(\pi^{-i},\sigma)}^{\XX\times\AA^i}$ and $\gamma'= \mu_{\eta,(\pi^{-i},\sigma')}^{\XX\times\AA^i}$.
Convexity of $\tilde{\bcal{O}}_\eta$  implies that 
\begin{equation}\label{eq-convexity-1}
\mu_{\eta,\hat{\pi}} =\alpha \mu_{\eta,(\pi^{-i},\sigma)}+(1-\alpha) \mu_{\eta,(\pi^{-i},\sigma')}
\end{equation}
for some $\hat{\pi}\in\tilde{\bcal{Y}}$. To get the result, let us show that $\mu^{{\XX\times\AA^{i}}}_{\eta,\hat{\pi}}=\mu^{{\XX\times\AA^{i}}}_{\eta,(\pi^{-i},\tilde{\sigma})}$
for some $\tilde{\sigma}\in\bcal{Y}^i$.
Observe that, by \eqref{eq-convexity-1} and Proposition \ref{prop-linear-equation-occupation}(ii),
\begin{align*}
\mu_{\eta,\hat{\pi}}  & =\big[\alpha  \mu_{\eta,(\pi^{-i},\sigma)}^{\XX\times\AA^{i}} + (1-\alpha) \mu_{\eta,(\pi^{-i},\sigma')}^{\XX\times\AA^{i}}\big]\otimes\pi^{-i}
= \mu^{\XX\times\AA^{i}}_{\eta,\hat{\pi}}  \otimes\pi^{-i}.
\end{align*}
Moreover, $\mu^{\XX\times\AA^{i}}_{\eta,\hat{\pi}}=\mu^{\XX}_{\eta,\hat{\pi}}\otimes\hat{\pi}^{\AA_{i}}$ and so, letting $\tilde\sigma=\hat\pi^{\AA^i}\in\bcal{Y}^i$ we obtain
$\mu_{\eta,\hat{\pi}}  =\mu_{\eta,\hat{\pi}} ^\XX\otimes(\pi^{-i},\tilde\sigma)$.
Since $\mu_{\eta,\hat{\pi}} ^\XX$ is equivalent to $\lambda$ on $\Delta^c$ (recall Proposition \ref{prop-linear-equation-occupation}(iv)), it follows that $\hat\pi$ and $(\pi^{-i},\tilde\sigma)$ coincide $\lambda$-a.s. on $\Delta^c$. In particular $\gamma_{\hat\pi}=\gamma_{(\pi^{-i},\tilde{\sigma})}$ and thus, using Lemma \ref{lemma-preliminary-1}(i),
$$\mu_{\eta,\hat{\pi}}=\mu_{\eta,\gamma_{\hat{\pi}}}=\mu_{\eta,\gamma_{(\pi^{-i},\tilde{\sigma})}}=\mu_{\eta,(\pi^{-i},\tilde\sigma)},$$ as we wanted to prove. This establishes convexity of $\mathcal{L}_{\eta}^i(\pi^{-i})$.
To prove compactness we will show that $\mathcal{L}_{\eta}^i(\pi^{-i})$ is closed. Suppose that $\gamma_n\rightarrow\gamma$ where $\gamma_n\in\mathcal{L}_{\eta}^{i}(\pi^{-i})$ and $\gamma\in\tilde{\bcal{O}}_\eta^i$. For each $n$ there is some $\sigma^i_n\in\bcal{Y}^i $ such that 
$\gamma_n=\mu_{\eta,(\pi^{-i},\sigma^i_n)}^{\XX\times\AA^i}$.
For some subsequence of $\{\sigma^i_n\}$, still denoted by $\{\sigma^i_n\}$, we have $\sigma^i_n\rightarrow \sigma_*^i$ for some $\sigma_*^i\in\bcal{Y}^i$.
By Corollary \ref{Convergence-mu}(iii), this shows that $\gamma_n\rightarrow \mu^{\XX\times\AA^i}_{\eta,(\pi^{-i},\sigma^i_*)}$, which indeed belongs to $\mathcal{L}_{\eta}^i(\pi^{-i})$. 
\hfill$\Box$\\[10pt]\indent
Given a player $1\le i\le N$ and a strategy profile $\pi^{-i}\in\bcal{Y}^{-i}$ for the remaining players,  let
$$
\mathcal{A}^i_{\eta,\rho^i}(\pi^{-i}) 
 = \Big\{ \mu_{\eta,(\pi^{-i},\sigma)}^{\XX\times\AA^{i}} : \sigma\in\bcal{Y}^i\hbox{ is such that } C^i(\eta,(\pi^{-i},\sigma))\ge\rho^i\Big\}\subseteq\mathcal{L}^i(\pi^{-i}).$$
 Thus, $\mathcal{A}^i_{\eta,\rho^i}(\pi^{-i})$ is the set of $(\XX\times\AA^i)$-marginals of the occupation measures of the Markov policies $\sigma\in\bcal{Y}^i$ of player $i$ such that the Markov profile $(\pi^{-i},\sigma)$ satisfies player $i$'s  constraint.

\begin{proposition}\label{prop-helps-lower}
Consider a player $1\le i\le N$ and a sequence $\{\pi^{-i}_n\}\subseteq\bcal{Y}^{-i}$ such that $\pi_n^{-i}\rightarrow\pi^{-i}$ in the product topology of $\bcal{Y}^{-i}$
for some $\pi\in\bcal{Y}$.
\begin{enumerate}
\item[(i).] If $C^i(\eta,\pi)\ge\rho^i$ then there exists a sequence $\{\gamma_n\}$ in $\bcal{M}^+(\XX\times\AA^i)$ such that $\gamma_n\rightarrow \mu_{\eta,\pi}^{\XX\times\AA^i}$ and for some $K\in\NN$, we have $\gamma_n\in\bcal{A}^i_{\eta,\rho^i}(\pi_n^{-i})$ for $n\ge K$.
\item[(ii).] If $\mathcal{G}(\eta,\rho)$ satisfies the Slater condition and $C^i(\eta,\pi)\ge\rho^i$ then there exists a sequence $\{\gamma_n\}$ in $\bcal{M}^+(\XX\times\AA^i)$ such that $\gamma_n\rightarrow \mu_{\eta,\pi}^{\XX\times\AA^i}$ and such that, for some $K\in\NN$, we have $\gamma_n\in\bcal{A}^i_{\eta_n,\rho^i}(\pi_n^{-i})$ for $n\ge K$.
\end{enumerate}
\end{proposition}
\textbf{Proof:} 
We will prove only item (ii) since item (i) can be easily obtained by using the same arguments.
From Corollary \ref{Convergence-mu}(iv) we have  
$\lim_{n\rightarrow\infty} C^i(\eta_{n},(\pi_n^{-i},\pi^i))= C^i(\eta,\pi)\ge\rho^{i}$.
So, there exist some sequence $\{\epsilon_{n}\}_{n\ge1}$ contained in $ [0,1)$ with $\epsilon_n\rightarrow0$ and some index $n_0$ for which
\begin{equation*}
C^i(\eta_{n},(\pi_n^{-i},\pi^i))\ge \rho^i-\epsilon_n\mathbf{1}
\quad\hbox{for all $n\ge n_0$}.
\end{equation*}
By the Slater condition,
 we can find  some $\bar\pi^i\in\bcal{Y}^i$ and $\delta>0$ such that 
$C^i(\eta,(\pi^{-i},\bar\pi^i))>\rho^{i}+\delta\mathbf{1}$.
Again from Corollary \ref{Convergence-mu}(iv),
there is some $n_1\ge n_0$ such that
\begin{equation*}
C^i(\eta_{n},(\pi_n^{-i},\bar\pi^i)) >\rho^{i}+\delta\mathbf{1}\quad\hbox{for all $n\ge n_1$}.
\end{equation*}
Observe  that for any $n\ge1$  both 
$\mu^{\XX\times\AA^i}_{\eta_{n},(\pi_n^{-i},\pi^i)}$ and $\mu^{\XX\times\AA^i}_{\eta_{n},(\pi_n^{-i},\bar\pi^i)}$ belong to $\mathcal{L}^i_{\eta_{n}}(\pi_n^{-i})$
which is a convex set by Proposition \ref{prop-L1-convex-compact}.
Hence, there exists some $\gamma_n\in\mathcal{L}^i_{\eta_n}(\pi^{-i}_n)$ and  $\sigma^i_n\in\bcal{Y}^i$ such that 
\begin{align}
\gamma_n=\mu^{\XX\times\AA^i}_{\eta_{n},(\pi^{-i}_{n} ,\sigma^i_n)} = (1-\sqrt{\epsilon_{n}}) \mu^{\XX\times\AA^i}_{\eta_{n},(\pi^{-i}_{n} ,\pi^i)}
+ \sqrt{\epsilon_{n}} \mu^{\XX\times\AA^i}_{\eta_{n},(\pi^{-i}_{n} ,\bar\pi^i)}
\label{eq-Occup-measure-def-pi1} 
\end{align}
with, as a consequence of Proposition \ref{prop-linear-equation-occupation}(ii),
\begin{equation*} 
\mu_{\eta_{n},(\pi^{-i}_{n} ,\sigma^i_n)}
=\mu^{\XX\times\AA^i}_{\eta_{n},(\pi^{-i}_{n} ,\sigma^i_n)}\otimes \pi^{-i}_n
 = (1-\sqrt{\epsilon_{n}}) \mu_{\eta_{n},(\pi^{-i}_{n} ,\pi^i)}
+ \sqrt{\epsilon_{n}} \mu_{\eta_{n},(\pi^{-i}_{n} ,\bar\pi^i)}
\end{equation*}
and so for any $n\ge n_1$
\begin{eqnarray*}
C^i(\eta_n,(\pi^{-i}_n,\sigma^i_n)) &=& (1-\sqrt{\epsilon_n}) C^i(\eta_n,(\pi^{-i}_n,\pi^i))+\sqrt{\epsilon_n} C^i(\eta_n,(\pi^{-i}_n,\bar\pi^i))\nonumber \\ \noalign{\smallskip}
&\ge&
\rho^i + \sqrt{\epsilon_{n}} \big[\delta -(1-\sqrt{\epsilon_{n}}) \sqrt{\epsilon_{n}} \big]\mathbf{1}.
\end{eqnarray*}
Therefore, there exists some $K\ge n_1$ such that $n\ge K$ implies $C^i(\eta_n,(\pi^{-i}_n,\sigma^i_n))\ge\rho^i$.
Since  $\gamma_n=\mu^{\XX\times\AA^i}_{\eta_{n},(\pi^{-i}_n, \sigma^i_n)}\in\mathcal{L}^i_{\eta_n}(\pi^{-i}_n)$,
this establishes precisely that $\gamma_n\in\mathcal{A}^i_{\eta_n,\rho^i}(\pi^{-i}_n)$ for all $n\ge K$. 
Finally, from \eqref{eq-Occup-measure-def-pi1} and Corollary \ref{Convergence-mu}(iii) we have that 
$\lim_{n\rightarrow\infty} \gamma_n=\lim_{n\rightarrow\infty} \mu^{\XX\times\AA^i}_{\eta_n,(\pi_n^{-i},\pi^i)}= \mu_{\eta,\pi}^{\XX\times\AA^i}$.
This completes the proof.
\hfill $\Box$

\begin{proposition}
\label{Continuity-A-1}
Given any $1\le i\le N$,  
the correspondence $\mathcal{A}^i_{\eta,\rho^i}:\bcal{Y}^{-i}\tto\tilde{\bcal{O}}_\eta^i$ defined by 
$\pi^{-i}\mapsto \mathcal{A}^i_{\eta,\rho^i}(\pi^{-i})$ is continuous with nonempty convex and compact values.
\end{proposition}
\textbf{Proof.} The Slater condition implies that $\mathcal{A}^i_{\eta,\rho^i}(\pi^{-i})$ is nonempty for any $\pi^{-i}\in\bcal{Y}^{-i}$.
To prove convexity, let $\gamma,\gamma'\in\mathcal{A}^i_{\eta,\rho^i}(\pi^{-i})$ and $0\le \alpha\le1$. Then, there exist $\sigma,\sigma'\in\bcal{Y}^i$ with
such that $\gamma=\mu^{\XX\times\AA^i}_{\eta,(\pi^{-i},\sigma)}$ and $\gamma'=\mu^{\XX\times\AA^i}_{\eta,(\pi^{-i},\sigma')}$
and satisfying $C^i(\eta,(\pi^{-i},\sigma))\ge\rho^i$ and  $C^i(\eta,(\pi^{-i},\sigma'))\ge\rho^i$.
By convexity of $\mathcal{L}^i_{\eta}(\pi^{-i})$ in Proposition \ref{prop-L1-convex-compact}, there exists some $\sigma^*\in\bcal{Y}^i$ such that
$\mu^{\XX\times\AA^i}_{\eta,(\pi^{-i},\sigma^*)} =\alpha \gamma+(1-\alpha)\gamma'$.
But then
\begin{eqnarray*}
\mu_{\eta,(\pi^{-i},\sigma^*)}&=& \mu^{\XX\times\AA^i}_{\eta,(\pi^{-i},\sigma^*)}\otimes \pi^{-i} = \alpha( \gamma\otimes\pi^{-i})+(1-\alpha)(\gamma'\otimes\pi^{-i})\\
&=&\alpha \mu_{\eta,(\pi^{-i},\sigma)} +(1-\alpha)\mu_{\eta,(\pi^{-i},\sigma')},
\end{eqnarray*}
and so by integration of the function $c^i(x,a)$ with respect to these occupation measures,
$$C^i(\eta,(\pi^{-i},\sigma^*))=\alpha  C^i(\eta,(\pi^{-i},\sigma))+(1-\alpha) C^i(\eta,(\pi^{-i},\sigma'))\ge\rho^i,$$
which establishes that $\alpha\gamma+(1-\alpha)\gamma'$ is in $\mathcal{A}^i_{\eta,\rho^i}(\pi^{-i})$.

The correspondence $\mathcal{A}^i_{\eta,\rho^i}$ takes values in the compact metric space $\tilde{\bcal{O}}^i_\eta$ and thus, by the Closed Graph Theorem in \cite[Theorem 17.11]{aliprantis06}, it is upper semicontinuous and compact-valued if and only if its graph is closed. 
Suppose that we have a convergent sequence $(\pi_n^{-i},\gamma_n)$ in the graph of $\mathcal{A}^i_{\eta,\rho^i} $ converging to some $(\pi^{-i},\gamma)\in\bcal{Y}^{-i}\times\tilde{\bcal{O}}_\eta^i$. We must show that $\gamma\in\mathcal{A}^i_{\eta,\rho^i}(\pi^{-i})$.
For each $n\ge1$ there exists $\sigma_n\in \bcal{Y}^{i}$ such that
\begin{align}
\gamma_{n} = \mu_{\eta,(\pi^{-i}_n,\sigma_n)}^{\XX\times\AA^{i}} \quad\hbox{and}\quad C^i(\eta,(\pi^{-i}_n,\sigma_n))\ge\rho^i.
\label{eq-tool-again}
\end{align}
For some subsequence $\{\sigma_{n'}\}$ of $\{\sigma_{n}\}$ we have $\sigma_{n'}\rightarrow\sigma$ for some $\sigma\in\bcal{Y}^i$ and so using Corollary \ref{Convergence-mu}(iii) and (iv)
$$ \gamma=\mu_{\eta,(\pi^{-i},\sigma)}^{\XX\times\AA^i}\quad\hbox{and}\quad C^i(\eta,(\pi^{-i},\sigma))\ge\rho^i,$$
as we wanted to prove. 
Lower semicontinuity of the correspondence follows from Proposition \ref{prop-helps-lower}(i) and  the sequential characterization of lower semicontinuity  in  \cite[Theorem 17.21]{aliprantis06}. 
\hfill $\Box$
\\[10pt]\indent
Given some player $1\le i\le N$ and a Markov profile $\pi^{-i}\in\bcal{Y}^{-i}$ for the remaining players, if player $i$ chooses the Markov policy $\sigma\in\bcal{Y}^i$ then his payoff is (recall Proposition \ref{prop-linear-equation-occupation}(ii))
$$R^i(\eta,(\pi^{-i},\sigma))=\int_{\XX\times\AA} r^id\mu_{\eta,(\pi^{-i},\sigma)}=
\int_{\XX\times\AA} r^i d(\mu_{\eta,(\pi^{-i},\sigma)}^{\XX\times\AA^i}\otimes\pi^{-i}),$$
and his goal is to maximize this payoff over all $\sigma\in\bcal{Y}^i$ such that $C^i(\eta,(\pi^{-i},\sigma))\ge\rho^i$ or, which is the same, 
maximize
$$\int_{\XX\times\AA} r^i d(\gamma\otimes\pi^{-i})$$
over all $\gamma$ belonging to $\mathcal{A}_{\eta,\rho^i}^i(\pi^{-i})$.  Based on this, we define the correspondence 
 $\mathcal{H}^i_{\eta,\rho^i}:\bcal{Y}^{-i}\tto \tilde{\bcal{O}}^i_{\eta}$ given by
\begin{equation}\label{eq-defin-H}
 \mathcal{H}^i_{\eta,\rho^i}(\pi^{-i})=\argmax_{\gamma\in\mathcal{A}^i_{\eta,\rho^i}(\pi^{-i})}\Big\{ \int_{\XX\times\AA} r^i d(\gamma\otimes\pi^{-i})\Big\}.
\end{equation}

\begin{proposition}
\label{Correspondence-H1}
For any $i\in\{1,\ldots,N\}$, 
the correspondence $\mathcal{H}^i_{\eta,\rho^i}:\bcal{Y}^{-i}\tto\tilde{\bcal{O}}_\eta^i$ is upper semicontinuous with nonempty compact and convex values.
\end{proposition}
\textbf{Proof.} On the graph of the correspondence $\mathcal{A}^i_{\eta,\rho^i}$ consider the function
$$f_\eta(\pi^{-i},\gamma)=\int_{\XX\times\AA} r^i d(\gamma\otimes\pi^{-i})$$
and let us prove that it is continuous. To this end, let $\pi^{-i}_n\rightarrow\pi^{-i}$ in $\bcal{Y}^{-i}$ and $\gamma_n\rightarrow\gamma$ in $\tilde{\bcal{O}}_\eta^i$ with $\gamma_n\in \mathcal{A}^i_{\eta,\rho^i}(\pi_n^{-i})$ and $\gamma\in \mathcal{A}^i_{\eta,\rho^i}(\pi^{-i})$. We must show that $f_\eta(\pi^{-i}_n,\gamma_n)\rightarrow f_\eta(\pi,\gamma)$, and we will prove that this limit holds through any convergence subsequence, which will be denoted by $\{n\}$ without loss of generality.
There exist $\sigma_n,\sigma\in\bcal{Y}^i$ such that 
$$\gamma_n=\mu_{\eta,(\pi_{n}^{-i},\sigma_n)}^{\XX\times\AA^i}\quad\hbox{and}\quad
\gamma=\mu_{\eta,(\pi^{-i},\sigma)}^{\XX\times\AA^i}$$
and we can also assume that $\sigma_n\rightarrow\sigma^*$ for some $\sigma^*\in\bcal{Y}^i$. Using Corollary \ref{Convergence-mu}(iv) we obtain that 
$$R^i(\eta,(\pi^{-i}_n,\sigma_n))\rightarrow R^i(\eta,(\pi^{-i},\sigma^*)).$$
On the other hand, by Corollary \ref{Convergence-mu}(iii),
$$\gamma_n\rightarrow \mu^{\XX\times\AA^i}_{\eta,(\pi^{-i},\sigma^*)}\quad\hbox{and so}\quad \gamma=\mu^{\XX\times\AA^i}_{\eta,(\pi^{-i},\sigma^*)}=
\mu^{\XX\times\AA^i}_{\eta,(\pi^{-i},\sigma)}.$$
This shows that 
$$\mu_{\eta,(\pi^{-i},\sigma)}=\mu_{\eta,(\pi^{-i},\sigma)}^{\XX\times\AA^i}\otimes \pi^{-i}=
\mu_{\eta,(\pi^{-i},\sigma^*)}^{\XX\times\AA^i}\otimes \pi^{-i}=
\mu_{\eta,(\pi^{-i},\sigma^*)}$$
and so $R^i(\eta,(\pi^{-i}_n,\sigma_n))\rightarrow R^i(\eta,(\pi^{-i},\sigma))$, which can be also written as
$f_\eta(\gamma_n,\pi^{-i}_n)\rightarrow f_\eta(\gamma,\pi^{-i})$.
Once we know that $f_\eta$ is continuous on the graph of $\mathcal{A}^i_{\eta,\rho^i}$, 
we can apply the  Berge Maximum Theorem  \cite[Theorem 17.31]{aliprantis06} and conclude that  the $\argmax$ correspondence   $\mathcal{H}^i_{\eta,\rho^i}$   is upper semicontinuous with nonempty compact values.

Finally, observe that the function $f_\eta$ is linear in $\gamma$ for fixed $\pi^{-i}\in\bcal{Y}^{-i}$
and so the set of maximizers $\mathcal{H}^i_{\eta,\rho^i}(\pi^{-i})$ is convex. 
\hfill $\Box$
\\[10pt]\indent
By considering the product of the correspondences $\mathcal{H}^i_{\eta,\rho^i}$ we obtain 
the following result, which easily follows from \cite[Theorem 17.28]{aliprantis06}.
\begin{corollary}
The correspondence  $\mathcal{H}_{\eta,\rho}:\bcal{Y}\tto \tilde{\bcal{O}}^1_{\eta}\times\ldots\times \tilde{\bcal{O}}^N_{\eta}$
defined by
$$\pi\mapsto \prod_{i=1}^N  \mathcal{H}_{\eta,\rho}^i(\pi^{-i}).$$
is upper semicontinuous  with nonempty compact and convex values.
\end{corollary}

\paragraph{Proof of Proposition \ref{th-main-0}:}
To get the result, let us show that the following results hold.
\begin{itemize} 
\item[(i)]  The  correspondence 
$$\mathcal{H}_{\eta,\rho}\circ \mathcal{J}_\eta: \tilde{\bcal{O}}_\eta^1\times\ldots \times\tilde{\bcal{O}}_\eta^N\tto  \tilde{\bcal{O}}_\eta^1\times\ldots \times\tilde{\bcal{O}}_\eta^N
$$
has a fixed point $(\gamma_*^1,\ldots,\gamma^N_*)$.
\item[(ii)] The Markov profile $\pi_*\in\bcal{Y}$ given by $\pi_*^i=\mathcal{J}^i_\eta(\gamma_*^i)$ for $1\le i\le N$ is a constrained equilibrium in the class of all strategy profiles $\mathbf{\Pi}$ of the players for the game model $\mathcal{G}(\eta,\rho)$.
\end{itemize}
Let us first proceed to the proof of item (i). By \cite[Theorem 17.23]{aliprantis06}, the composition  $\mathcal{H}_{\eta,\rho}\circ\mathcal{J}_{\eta}$ is an upper semicontinuous correspondence as it is the composition of a continuous function $\mathcal{J}_{\eta}$ and an upper semicontinuous correspondence $\mathcal{H}_{\eta,\rho}$. Besides, it has nonempty compact and convex values.  By the Closed Graph Theorem \cite[Theorem 17.11]{aliprantis06}, it is also a closed  correspondence.

Since  $\tilde{\bcal{O}}^1_{\eta}\times\ldots\times\tilde{\bcal{O}}^N_{\eta}$ is a nonempty compact convex subset  of the locally convex Hausdorff space
$\bcal{M}(\XX\times\AA^{1})\times\ldots\times\bcal{M}(\XX\times\AA^N)$ ---Proposition 2.2 in \cite{dufour22}--- we can use the
 Kakutani-Fan-Glicksberg fixed point theorem \cite[Corollary 17.55]{aliprantis06}  to get the existence of a fixed point for  the correspondence $\mathcal{H}_{\eta,\rho}\circ\mathcal{J}_{\eta}$.\\[5pt]
\noindent
Let us now proceed to the proof of item (ii). If $(\gamma_*^1,\ldots,\gamma_*^N)$ is a fixed point of $\mathcal{H}_{\eta,\rho}\circ\mathcal{J}_{\eta}$, consider the Markov policies $\pi_*^i=\mathcal{J}^i_\eta(\gamma^i_*)\in\bcal{Y}^i$ for $1\le i\le N$ and let $\pi_*=(\pi^1_*,\ldots,\pi^N_*)\in\bcal{Y}$. Since $\gamma^i_*\in \mathcal{L}^i_{\eta}(\pi_*^{-i})$ we have that for some $\sigma\in\bcal{Y}^i$ it is 
$\gamma_*^i=\mu^{\XX\times\AA^i}_{\eta,(\pi_*^{-i},\sigma)}$. But now using Lemma \ref{lemma-preliminary-1}(iii) we also have
$$\gamma_*^i=\mu^{\XX\times\AA^i}_{\eta,(\pi_*^{-i},\sigma)}=\mu^{\XX\times\AA^i}_{\eta,(\pi_*^{-i},\mathcal{J}_\eta^i(\gamma^i_*))}=\mu^{\XX\times\AA^i}_{\eta,\pi^*}.$$
In particular, for each $1\le i\le N$ we have $\gamma^i_*\otimes\pi_*^{-i}=\mu_{\eta,\pi_*}$.
Moreover, since $\gamma_*^i\in\mathcal{A}^i_{\eta,\rho^i}(\pi^{-i}_*)$ it follows that $C^i(\eta,\pi_*)\ge\rho^i$. We conclude that the Markov profile $\pi_*\in\bcal{Y}$ satisfies the constraints of all the players.

If player $i$ varies his policy from $\pi_*^i\in\bcal{Y}^i$ to some $\pi^i\in\mathbf{\Pi}^i$ which satisfies his own constraint (i.e. $C^i(\eta,(\pi_*^{-i},\pi^i))\ge\rho^i$)
we can use the result in Proposition \ref{prop-linear-equation-occupation}(iii) to derive the existence of some $\sigma\in\bcal{Y}^i$ such that 
$$\mu_{\eta,(\pi_*^{-i},\pi^i)}=\mu_{\eta,(\pi_*^{-i},\sigma)}$$
with, again, $C^i(\eta,(\pi_*^{-i},\sigma))\ge\rho^i$. This implies that $\mu^{\XX\times\AA^i}_{\eta,(\pi_*^{-i},\sigma)}\in\mathcal{A}^i_{\eta,\rho^i}(\pi_*^{-i})$ and thus
$$\int_{\XX\times\AA} r^i d( \mu^{\XX\times\AA^i}_{\eta,(\pi_*^{-i},\sigma)}\otimes\pi_*^{-i}) \le
\int_{\XX\times\AA} r^i d (\gamma_*^i\otimes\pi_*^{-i}) $$ or, equivalently,
$$R^i(\eta,(\pi^{-i}_*,\pi))=R^i(\eta,(\pi_*^{-i},\sigma))\le R^i(\eta,\pi_*).$$
This completes the proof. 
\hfill$\Box$

\subsection{Proof of Theorem \ref{th-main}}
\label{Proof-Main2}
Clearly, for all $n\ge1$ we have $\lambda\ll \eta_n$.  We  also have that $\{\eta_n\}_{n\in\NN}$ converges to $\eta$ in total variation and that the corresponding density functions with respect to $\lambda$
converge strongly (or, in norm) in  $L^1(\XX,\bfrak{X},\lambda)$: $\ds \left\|{d\eta_n}/{d\lambda}-{d\eta}/{d\lambda}\right\|_1\rightarrow0$.
Since the constraint function $c_i$ is bounded by $\mathbf{r}$, we have $|C^i(\eta_n,\pi)-C^i(\eta,\pi)|\le \mathbf{r}/(n+1)$.
Also, for the constraint constants $\rho_n^i=\rho^i-\frac{\mathbf{r}}{n+1}\mathbf{1}$ with $1\le i\le N$,
the game model $\mathcal{G}(\eta_n,\rho_n)$ satisfies the Slater condition in Definition \ref{def-slater}. 
Under assumptions \ref{Assump-absorbing-uniformly} and \ref{Assump-absorbing-lambda-uniformly}, we obtain that  the game model
$\mathcal{G}(\eta_{n},\rho_{n})$ is uniformly absorbing to~$\Delta$ by using item (ii) of Proposition \ref{prop-preliminary}.
We can conclude that the game model $\mathcal{G}(\eta_{n},\rho_{n})$ satisfies Assumptions~\ref{Assump-2} and~\ref{Assumption-B}
and so,
Proposition \ref{th-main-0} yields the existence of a constrained Nash equilibrium $\hat{\pi}_{n}\in\bcal{Y}$ for the game model $\mathcal{G}(\eta_n,\rho_n)$ with $n\ge1$.
This means that 
\begin{equation}\label{eq-approx-game-1}
C^i(\eta_n,\hat{\pi}_{n})\ge\rho^i_n\quad\hbox{for $1\le i\le N$}
\end{equation}
and that, for any $1\le i\le N$ and $\pi^i\in\mathbf{\Pi}^i$,
$$C^i(\eta_n,(\hat{\pi}_{n}^{-i},\pi^i))\ge\rho^i_n\ \Rightarrow\ R^i(\eta_n,(\hat{\pi}_{n}^{-i},\pi^i))\le R^i(\eta_n,\hat{\pi}_{n}).$$

Without loss of generality, we assume that the sequence of so-defined equilibria converges to some $\hat{\pi}\in\bcal{Y}$, that is, for each $1\le i\le N$ we have
$\hat{\pi}_{n}^i\rightarrow\hat{\pi}^i$ in $\bcal{Y}^i$ as $n\rightarrow\infty$.
We want to show that $\hat{\pi}$ is a constrained equilibrium for the game model $\mathcal{G}(\eta,\rho)$.
To see this, note first that we can take the limit in \eqref{eq-approx-game-1} to obtain that $C^i(\eta,\hat{\pi})\ge\rho^i$ for every $1\le i\le N$
by using Corollary \ref{Convergence-mu}(iv).
Secondly, fix $i\in\{1,\ldots,N\}$ and choose any $\pi^i\in\mathbf{\Pi}^i$ such that $C^i(\eta,(\hat{\pi}^{-i},\pi^i))\ge\rho^i$.
By Proposition \ref{prop-linear-equation-occupation}(iii) it follows that there is some $\sigma\in\bcal{Y}^i$ such that $(\hat{\pi}^{-i},\pi^i)$ and $(\hat{\pi}^{-i},\sigma)\in\bcal{Y}$
yield the same payoffs $C^i$ and $R^i$. Hence we have  $C^i(\eta,(\hat{\pi}^{-i},\sigma))\ge\rho^i$ and we must show that
$$ R^i(\eta,(\hat{\pi}^{-i},\sigma))\le R^i(\eta,\hat{\pi}).$$
We will use Proposition \ref{prop-helps-lower}(ii) for the Markov profile $(\pi_*^{-i},\sigma)\in\bcal{Y}$ and the sequence $\{\hat{\pi}_{n}^{-i}\}$
to derive the existence of a sequence
$\gamma_n\rightarrow \mu^{\XX\times\AA^i}_{\eta,(\hat{\pi}^{-i},\sigma)}$ such that $\gamma_n\in\mathcal{A}^i_{\eta_n,\rho^i}(\hat{\pi}_{n}^{-i})$  for large enough $n\ge K$.  So, for such  $n\ge K$, let $\sigma_n\in\bcal{Y}^i$  be such that $\gamma_n=\mu^{\XX\times\AA^i}_{\eta_n,(\hat{\pi}_{n}^{-i},\sigma_n)}$ which satisfies
$C^i(\eta_n,(\hat{\pi}_{n}^{-i},\sigma_n))\ge\rho^i\ge \rho^i_n$.
This implies that for any $n\ge K$ we have
$$R^i(\eta_n,(\hat{\pi}_{n}^{-i},\sigma_n))\le R^i(\eta_n,\hat{\pi}_{n}).$$
There exists some $\bar\sigma\in\bcal{Y}^i$ and a subsequence of $\{\sigma_n\}$ (still denoted by $\{\sigma_n\}$) satisfying $\sigma_n\rightarrow\bar\sigma$ in $\bcal{Y}^i$ and then taking the limit we have
$$R^i(\eta,(\hat{\pi}^{-i},\bar\sigma))\le R^i(\eta,\hat{\pi})$$
by using Corollary \ref{Convergence-mu}(iv).
But then item (iii) of Corollary \ref{Convergence-mu} implies that 
$$\gamma_n=\mu^{\XX\times\AA^i}_{\eta_n,(\hat{\pi}_{n}^{-i},\sigma_n)}\rightarrow
\mu^{\XX\times\AA^i}_{\eta,(\hat{\pi}^{-i},\bar\sigma)}= \mu^{\XX\times\AA^i}_{\eta,(\hat{\pi}^{-i},\sigma)}$$
so that 
$$\mu_{\eta,(\hat{\pi}^{-i},\bar\sigma)}= \mu^{\XX\times\AA^i}_{\eta,(\hat{\pi}^{-i},\bar\sigma)}\otimes \hat{\pi}^{-i}
=\mu^{\XX\times\AA^i}_{\eta,(\hat{\pi}^{-i},\sigma)}\otimes \hat{\pi}^{-i}= \mu_{\eta,(\hat{\pi}^{-i},\sigma)}$$ and, hence, 
$R^i(\eta,(\hat{\pi}^{-i},\bar\sigma))= R^i(\eta,(\hat{\pi}^{-i},\sigma))$
and $R^i(\eta,(\hat{\pi}^{-i},\sigma))\le R^i(\eta,\hat{\pi})$ follows.

\subsection{Proof of  Corollary \ref{Cond-sufficient-Assumption-B}}\label{sec53}
(i).  
To check this result, we must show that the convergences 
$\pi^i_n\rightarrow \pi^i$ in $\bcal{Y}^i$ for each $1\le i\le N$ imply that 
$$\pi_n(da|x)=\pi^1_{n}(da^1|x)\times\cdots\times\pi^N_{n}(da^N|x)\rightarrow \pi^1(da^1|x)\times\cdots\times\pi^N(da^N|x)=\pi(da|x)\quad\hbox{in $\tilde{\bcal{Y}}$}.$$ 
To avoid trivial cases, suppose that $\lambda\{x\}>0$ for every $x\in\XX$. Then, $\pi^i_n\rightarrow\pi^i$ means that $\pi^i_n(da|x)$ converges  in the weak topology of $\bcal{P}(\AA^i)$ to $\pi^i(da|x)$ for any $x\in\XX$. By \cite[Theorem 2.8]{billingsley} it follows that 
$$\pi^1_{n}(da^1|x)\times\cdots\times\pi^N_{n}(da^N|x)\rightarrow\pi^1(da^1|x)\times\cdots\times\pi^N(da^N|x)$$ 
in the weak topology of $\bcal{P}(\AA)$  for any $x\in\XX$. 
Given arbitrary $f\in\car(\XX\times\AA)$ bounded by a function  $F\in L^1(\XX,\bfrak{X},\lambda)$, that is, with $\sum_x F(x)\lambda\{x\}<\infty$, from the dominated convergence theorem we obtain that
$$\sum_{x\in\XX} \int_\AA f(x,a)\pi_n(da|x)\lambda\{x\}\rightarrow \sum_{x\in\XX} \int_\AA f(x,a)\pi(da|x)\lambda\{x\}.
$$
which shows that, indeed, $\pi_n\rightarrow\pi$ in $\tilde{\bcal{Y}}$. As a direct consequence of Lemma \ref{lem-trace-2}, we conclude that the continuity properties in Assumption \ref{Assumption-B} are satisfied.

Note that this proof establishes, in fact, that the trace topology of $\tilde{\bcal{Y}}$ on $\bcal{Y}$ coincides with the product topology of $\bcal{Y}=\bcal{Y}^1\times\ldots\times\bcal{Y}^N$. Such a result is known as a \textit{fiber product lemma}.
\par\noindent
(ii). Under the additive reward condition, the continuity of $\pi\mapsto r^i_\pi$ and $\pi\mapsto c^{i,j}_\pi$ is trivial since those functions turn out to be the sum of continuous functions. Regarding the additive transition property, observe that the density function
$$(y,x,a^1,\ldots,a^N)\mapsto \sum_{l=1}^N q^l(y,x,a^l)$$
satisfies the conditions in Assumption \ref{Assumption-transition-Q}. Checking the continuity of $\pi\mapsto Q_\pi v$ on $\bcal{Y}$ is again straightforward by using the additive property of the density function.
\hfill$\Box$

\end{document}